# GENERALIZED FUNCTIONAL LINEAR MODELS[1]

### By Hans-Georg Müller and Ulrich Stadtmüller

*University of California, Davis and Universität Ulm*

We propose a generalized functional linear regression model for a regression situation where the response variable is a scalar and the predictor is a random function. A linear predictor is obtained by forming the scalar product of the predictor function with a smooth parameter function, and the expected value of the response is related to this linear predictor via a link function. If, in addition, a variance function is specified, this leads to a functional estimating equation which corresponds to maximizing a functional quasi-likelihood. This general approach includes the special cases of the functional linear model, as well as functional Poisson regression and functional binomial regression. The latter leads to procedures for classification and discrimination of stochastic processes and functional data. We also consider the situation where the link and variance functions are unknown and are estimated nonparametrically from the data, using a semiparametric quasi-likelihood procedure.

An essential step in our proposal is dimension reduction by approximating the predictor processes with a truncated Karhunen–Loève expansion. We develop asymptotic inference for the proposed class of generalized regression models. In the proposed asymptotic approach, the truncation parameter increases with sample size, and a martingale central limit theorem is applied to establish the resulting increasing dimension asymptotics. We establish asymptotic normality for a properly scaled distance between estimated and true functions that corresponds to a suitable $L^2$ metric and is defined through a generalized covariance operator. As a consequence, we obtain asymptotic tests and simultaneous confidence bands for the parameter function that determines the model.

The proposed estimation, inference and classification procedures and variants with unknown link and variance functions are investi-

Received September 2001; revised March 2004.

[1]Supported in part by NSF Grants DMS-99-71602 and DMS-02-04869.

gated in a simulation study. We find that the practical selection of the number of components works well with the AIC criterion, and this finding is supported by theoretical considerations. We include an application to the classification of medflies regarding their remaining longevity status, based on the observed initial egg-laying curve for each of 534 female medflies.

**1. Introduction.** Many studies involve tightly spaced repeated measurements on the same individuals or direct recordings of a sample of curves [Brumback and Rice (1998) and Staniswalis and Lee (1998)]. If longitudinal measurements are made on a suitably dense grid, such data can often be regarded as a sample of curves or as functional data. Examples can be found in studies on longevity and reproduction, where typical subjects are fruit flies [Müller et al. (2001)] or nematodes [Wang, Müller, Capra and Carey (1994)].

Our procedures are motivated by a study where the goal is to find out whether there is information in the egg-laying curve observed for the first 30 days of life for female medflies, regarding whether the fly is going to be long-lived or short-lived. Discrimination and classification of curve data is of wide interest, from engineering [Hall, Poskitt and Presnell (2001)], and astronomy [Hall, Reimann and Rice (2000)] to DNA expression arrays with repeated measurements, where dynamic classification of genes is of interest [Alter, Brown and Botstein (2000)]. For multivariate predictors with fixed dimension, such discrimination tasks are often addressed by fitting binomial regression models using quasi-likelihood based estimating equations.

Given the importance of discrimination problems for curve data, it is clearly of interest to extend notions such as logistic, binomial or Poisson regression to the case of a functional predictor, which may be often viewed as a random predictor process. More generally, there is a need for new models and procedures allowing one to regress univariate responses of various types on a predictor process. The extension from the classical situation with a finite-dimensional predictor vector to the case of an infinite-dimensional predictor process involves a distinctly different and more complicated technology. One characteristic feature is that the asymptotic analysis involves increasing dimension asymptotics, where one considers a sequence of increasingly larger models.

The functional linear regression model with functional or continuous response has been the focus of various investigations [see Ramsay and Silverman (1997), Faraway (1997), Cardot, Ferraty and Sarda (1999) and Fan and Zhang (2000)]. An applied version of a generalized linear model with functional predictors has been investigated by James (2002). We assume here that the dependent variable is univariate and continuous or discrete, for example, of binomial or Poisson type, and that the predictor is a random function.

The main idea is to employ a Karhunen–Loève or other orthogonal expansion of the random predictor function [see, e.g., Ash and Gardner (1975) and Castro, Lawton and Sylvestre (1986)], with the aim to reduce the dimension to the first few components of such an expansion. The expansion is therefore truncated at a finite number of terms which increases asymptotically.

Once the dimension is reduced to a finite number of components, the expansion coefficients of the predictor process determine a finite-dimensional vector of random variables. We can then apply the machinery of generalized linear or quasi-likelihood models [Wedderburn (1974)], essentially solving an estimating or generalized score equation. The resulting regression coefficients obtained for the linear predictor in such a model then provide us with an estimate of the parameter function of the generalized functional regression model. This parameter function replaces the parameter vector of the ordinary finite-dimensional generalized linear model. We derive an asymptotic limit result (Theorem 4.1) for the deviation between estimated and true parameter function for increasing dimension asymptotics, referring to a situation where the number of components in the model increases with sample size.

Asymptotic tests for the regression effect and simultaneous confidence bands are obtained as corollaries of this main result. We include an extension to the case of a semiparametric quasi-likelihood regression (SPQR) model in which link and variance functions are unknown and are estimated from the data, extending previous approaches of Chiou and Müller (1998, 1999), and also provide an analysis of the AIC criterion for order selection.

The paper is organized as follows: The basics of the proposed generalized functional linear model and some preliminary considerations can be found in Section 2. The underlying ideas of estimation and statistical analysis within the generalized functional linear model will be discussed in Section 3. The main results and their ramifications are described in Section 4, preceded by a discussion of the appropriate metric in which to formulate the asymptotic result, which is found to be tied to the link and variance functions used for the generalized functional linear model. Simulation results are reported in Section 5. An illustrative example for the special case of binomial functional regression with the goal to discriminate between short- and long-lived medflies is provided in Section 6. This is followed by the main proofs in Section 7. Proofs of auxiliary results are in the Appendix.

**2. The generalized functional linear model.** The data we observe for the $i$th subject or experimental unit are $(\{X_i(t), t \in \mathcal{T}\}, Y_i)$, $i = 1, \ldots, n$. We assume that these data form an i.i.d. sample. The predictor variable $X(t)$, $t \in \mathcal{T}$, is a random curve which is observed per subject or experimental unit and corresponds to a square integrable stochastic process on a real interval $\mathcal{T}$. The dependent variable $Y$ is a real-valued random variable which may be continuous or discrete. For example, in the important special case of a binomial functional regression, one would have $Y \in \{0, 1\}$.

Assume that a link function $g(\cdot)$ is given which is a monotone and twice continuously differentiable function with bounded derivatives and is thus invertible. Furthermore, we have a variance function $\sigma^2(\cdot)$ which is defined on the range of the link function and is strictly positive. The generalized functional linear model or functional quasi-likelihood model is determined by a parameter function $\beta(\cdot)$, which is assumed to be square integrable on its domain $\mathcal{T}$, in addition to the link function $g(\cdot)$ and the variance function $\sigma^2(\cdot)$.

Given a real measure $dw$ on $\mathcal{T}$, define linear predictors

$$\eta = \alpha + \int \beta(t) X(t) \, dw(t)$$

and conditional means $\mu = g(\eta)$, where $\mathbf{E}(Y|X(t), t \in \mathcal{T}) = \mu$ and $\mathbf{Var}(Y|X(t), t \in \mathcal{T}) = \sigma^2(\mu) = \tilde{\sigma}^2(\eta)$ for a function $\tilde{\sigma}^2(\eta) = \sigma^2(g(\eta))$. In a generalized functional linear model the distribution of $Y$ would be specified within the exponential family. For the following (except where explicitly noted), it will be sufficient to consider the functional quasi-likelihood model

$$(1) \qquad Y_i = g\left(\alpha + \int \beta(t) X_i(t) \, dw(t)\right) + e_i, \qquad i = 1, \ldots, n,$$

where

$$\mathbf{E}(e|X(t), t \in \mathcal{T}) = 0,$$
$$\mathbf{Var}(e|X(t), t \in \mathcal{T}) = \sigma^2(\mu) = \tilde{\sigma}^2(\eta).$$

Note that $\alpha$ is a constant, and the inclusion of an intercept allows us to require $\mathbf{E}(X(t)) = 0$ for all $t$.

The errors $e_i$ are i.i.d. and we use integration w.r.t. the measure $dw(t)$ to allow for nonnegative weight functions $v(\cdot)$ such that $v(t) > 0$ for $t \in \mathcal{T}$, $v(t) = 0$ for $t \notin \mathcal{T}$ and $dw(t) = v(t) \, dt$; the default choice will be $v(t) = \mathbb{1}_{\{t \in \mathcal{T}\}}$. Nonconstant weight functions might be of interest when the observed predictor processes are function estimates which may exhibit increased variability in some regions, for example, toward the boundaries.

The parameter function $\beta(\cdot)$ is a quantity of central interest in the statistical analysis and replaces the vector of slopes in a generalized linear model or estimating equation based model. Setting $\sigma^2 = \mathbf{E}\{\tilde{\sigma}^2(\eta)\}$, we then find

$$\mathbf{Var}(e) = \mathbf{Var}\{\mathbf{E}(e|X(t), t \in \mathcal{T}\} + \mathbf{E}\{\mathbf{Var}(e|X(t), t \in \mathcal{T}\}$$
$$= \mathbf{E}\{\tilde{\sigma}^2(\eta)\} = \sigma^2,$$

as well as $\mathbf{E}(e) = 0$.

Let $\rho_j$, $j = 1, 2, \ldots$, be an orthonormal basis of the function space $L^2(dw)$, that is, $\int_{\mathcal{T}} \rho_j(t)\rho_k(t)\,dw(t) = \delta_{jk}$. Then the predictor process $X(t)$ and the parameter function $\beta(t)$ can be expanded into

$$X(t) = \sum_{j=1}^{\infty} \varepsilon_j \rho_j(t), \qquad \beta(t) = \sum_{j=1}^{\infty} \beta_j \rho_j(t)$$

[in the $L^2(dw)$ sense] with r.v.'s $\varepsilon_j$ and coefficients $\beta_j$, given by $\varepsilon_j = \int X(t) \times \rho_j(t)\,dw(t)$ and $\beta_j = \int \beta(t)\rho_j(t)\,dw(t)$, respectively. We note that $\mathbf{E}(\varepsilon_j) = 0$ and $\sum \beta_j^2 < \infty$. Writing $\sigma_j^2 = \mathbf{E}(\varepsilon_j^2)$, we find $\sum \sigma_j^2 = \int \mathbf{E}(X^2(t))\,dw(t) < \infty$.

From the orthonormality of the base functions $\rho_j$, it follows immediately that

$$\int \beta(t)X(t)\,dw(t) = \sum_{j=1}^{\infty} \beta_j \varepsilon_j.$$

It will be convenient to work with standardized errors

$$e' = e\sigma(\mu) = e\tilde{\sigma}(\eta),$$

for which $\mathbf{E}(e'|X) = 0$, $\mathbf{E}(e') = 0$, $\mathbf{E}(e'^2) = 1$. We assume that $\mathbf{E}(e'^4) = \mu_4 < \infty$ and note that in model (1), the distribution of the errors $e_i$ does not need to be specified and, in particular, does not need to be a member of the exponential family. In this regard, model (1) is less an extension of the classical generalized linear model [McCullagh and Nelder (1989)] than an extension of the quasi-likelihood approach of Wedderburn (1974). We address the difficulty caused by the infinite dimensionality of the predictors by approximating model (1) with a series of models where the number of predictors is truncated at $p = p_n$ and the dimension $p_n$ increases asymptotically as $n \to \infty$.

A heuristic motivation for this truncation strategy is as follows: Setting

$$U_p = \alpha + \sum_{j=1}^{p} \beta_j \varepsilon_j, \qquad V_p = \sum_{j=p+1}^{\infty} \beta_j \varepsilon_j,$$

we find $\mathbf{E}(Y|X(t), t \in \mathcal{T}) = g(\alpha + \sum_{j=1}^{\infty} \beta_j \varepsilon_j) = g(U_p + V_p)$. Conditioning on the first $p$ components and writing $F_{V_p|U_p}$ for the conditional distribution function leads to a truncated link function $g_p$,

$$\mathbf{E}(Y|U_p) = g_p(U_p) = \mathbf{E}[g(U_p + V_p)|U_p] = \int g(U_p + s)\,dF_{V_p|U_p}(s).$$

For the approximation of the full model by the truncated link function, we note that the boundedness of $g'$, $|g'(\cdot)|^2 \leq c$, implies that

$$\left\{ \int [g(U_p + V_p) - g(U_p + s)]\,dF_{V_p|U_p}(s) \right\}^2$$

$$\leq \int g'(\xi)^2 (V_p - s)^2 \, dF_{V_p|U_p}(s)$$

$$\leq 2c \int (V_p^2 + s^2) \, dF_{V_p|U_p}(s)$$

and, therefore,

$$\begin{aligned}
(2) \qquad &\mathbf{E}((g(U_p + V_p) - g_p(U_p))^2) \\
&= \mathbf{E}\left( \int [g(U_p + V_p) - g(U_p + s)] \, dF_{V_p|U_p}(s) \right)^2 \\
&\leq 2c\mathbf{E}(V_p^2 + \mathbf{E}(V_p^2|U_p)) = 4c\mathbf{E}(V_p^2) \\
&\leq 4c \sum_{j=p+1}^{\infty} \beta_j^2 \sum_{j=p+1}^{\infty} \sigma_j^2.
\end{aligned}$$

The approximation error of the truncated model is seen to be directly tied to $\mathbf{Var}(V_p)$ and is controlled by the sequence $\sigma_j^2 = \mathbf{Var}(\varepsilon_j)$, $j = 1, 2, \dots$, which for the special case of an eigenbase corresponds to a sequence of eigenvalues.

Setting $\varepsilon_j^{(i)} = \int X_i(t)\rho_j(t) \, dw(t)$, the full model with standardized errors $e_i'$ is

$$Y_i = g\left(\alpha + \sum_{j=1}^{\infty} \beta_j \varepsilon_j^{(i)}\right) + e_i'\tilde{\sigma}\left(\alpha + \sum_{j=1}^{\infty} \beta_j \varepsilon_j^{(i)}\right), \qquad i = 1, \dots, n.$$

With truncated linear predictors $\eta$ and means $\mu$,

$$\eta_i = \alpha + \sum_{j=1}^{p} \beta_j \varepsilon_j^{(i)}, \qquad \mu_i = g(\eta_i),$$

the $p$-truncated model becomes

$$(3) \quad Y_i^{(p)} = g_p\left(\alpha + \sum_{j=1}^{p} \beta_j \varepsilon_j^{(i)}\right) + e_i'\tilde{\sigma}_p\left(\alpha + \sum_{j=1}^{p} \beta_j \varepsilon_j^{(i)}\right), \qquad i = 1, \dots, n,$$

where $\tilde{\sigma}_p$ is defined analogously to $g_p$. Note that $g(U_p) - g_p(U_p)$ and, analogously, $\tilde{\sigma}(U_p) - \tilde{\sigma}_p(U_p)$ are bounded by the error (2). Since it will be assumed that this error vanishes asymptotically, as $p \to \infty$, we may instead of (3) work with the approximating sequence of models

$$(4) \quad Y_i^{(p)} = g\left(\alpha + \sum_{j=1}^{p} \beta_j \varepsilon_j^{(i)}\right) + e_i'\tilde{\sigma}\left(\alpha + \sum_{j=1}^{p} \beta_j \varepsilon_j^{(i)}\right), \qquad i = 1, \dots, n,$$

in which the functions $g$ and $\tilde{\sigma}$ are fixed. We note that the random variables $Y_i^{(p)}$ and $e_i'$, $i = 1, \dots, n$, form triangular arrays, $Y_{i,n}^{(p_n)}$ and $e_{i,n}'$, $i = 1, \dots, n$, with changing distribution as $n$ changes; for simplicity, we suppress the indices $n$.

Inference will be developed for the sequence of $p$-truncated models (4) with asymptotic results for $p \to \infty$. The practical choice of $p$ in finite sample situations will be discussed in Section 5. We also develop a version where the link function $g$ is estimated from the data, given $p$. The practical implementation of this semiparametric quasi-likelihood regression (SPQR) version adapts to the changing link functions $g_p$ of the approximating sequence (3).

**3. Estimation in the generalized functional linear model.** One central aim is estimation and inference for the parameter function $\beta(\cdot)$. Inference for $\beta(\cdot)$ is of interest for constructing confidence regions and testing whether the predictor function has any influence on the outcome, in analogy to the test for regression effect in a classical regression model. The orthonormal basis $\{\rho_j, j = 1, 2, \ldots\}$ is commonly chosen as the Fourier basis or the basis formed by the eigenfunctions of the covariance operator. The eigenfunctions can be estimated from the data as described in Rice and Silverman (1991) or Capra and Müller (1997). Whenever estimation and inference for the intercept $\alpha$ is to be included, we change the summation range for the linear predictors $\eta_i$ on the right-hand side of the $p$-truncated model (3) to $\sum_0^p$ from $\sum_1^p$, setting $\varepsilon_0^{(i)} = 1$ and $\beta_0 = \alpha$. In the following, inclusion of $\alpha$ into the parameter vector will be the default.

Fixing $p$ for the moment, we are in the situation of the usual estimating equation approach and can estimate the unknown parameter vector $\beta^T = (\beta_0, \ldots, \beta_p)$ by solving the estimating or score equation

$$(5) \qquad U(\beta) = 0.$$

Setting $\varepsilon^{(i)T} = (\varepsilon_0^{(i)}, \ldots, \varepsilon_p^{(i)})$, $\eta_i = \sum_{j=0}^p \beta_j \varepsilon_j^{(i)}$, $\mu_i = g(\eta_i)$, $i = 1, \ldots, n$, the vector-valued score function is defined by

$$(6) \qquad U(\beta) = \sum_{i=1}^n (Y_i - \mu_i) g'(\eta_i) \varepsilon^{(i)} / \sigma^2(\mu_i).$$

The solutions of the score equation (5) will be denoted by

$$(7) \qquad \hat{\beta}^T = (\hat{\beta}_0, \ldots, \hat{\beta}_p); \qquad \hat{\alpha} = \hat{\beta}_0.$$

Relevant matrices which play a well-known role in solving the estimating equation (5) are

$$D = D_{n,p} = (g'(\eta_i) \varepsilon_k^{(i)} / \sigma(\mu_i))_{1 \le i \le n, 0 \le k \le p},$$

$$V = V_{n,p} = \operatorname{diag}(\sigma^2(\mu_1), \ldots, \sigma^2(\mu_n))_{1 \le i \le n},$$

and with generic copies $\eta, \varepsilon, \mu$ of $\eta_i, \varepsilon^{(i)}, \mu_i$, respectively,

$$\Gamma = \Gamma_p = (\gamma_{kl})_{0 \leq k,l \leq p}, \qquad \gamma_{kl} = \mathbf{E}\left(\frac{g''^2(\eta)}{\sigma^2(\mu)} \varepsilon_k \varepsilon_l\right),$$

(8)

$$\Xi = \Gamma^{-1} = (\xi_{kl})_{0 \leq k,l \leq p}.$$

We note that $\Gamma = \frac{1}{n} E(D^T D)$ is a symmetric and positive definite matrix and that the inverse matrix $\Xi$ exists. Otherwise, one would arrive at the contradiction $\mathbf{E}((\sum_{k=0}^p \alpha_k \varepsilon_k g'(\eta)/\sigma(\mu))^2) = 0$ for nonzero constants $\alpha_0, \ldots, \alpha_p$.

With vectors $Y^T = (Y_1, \ldots, Y_n)$, $\mu^T = (\mu_1, \ldots, \mu_n)$, the estimating equation $U(\beta) = 0$ can be rewritten as

$$D^T V^{-1/2}(Y - \mu) = 0.$$

This equation is usually solved iteratively by the method of iterated weighted least squares. Under our basic assumptions, as $\frac{1}{n}\mathbf{E}(D^T D) = \Gamma_p$ is a fixed positive definite matrix for each $p$, the existence of a unique solution for each fixed $p$ is assured asymptotically.

In the above developments we have assumed that both the link function $g(\cdot)$ and the variance function $\sigma^2(\cdot)$ are known. Situations where the link and variance functions are unknown are common, and we can extend our methods to cover the general case where these functions are smooth, which for fixed $p$ corresponds to the semiparametric quasi-likelihood regression (SPQR) models considered in Chiou and Müller ([1998](#), [1999](#)). In the implementation of SPQR one alternates nonparametric (smoothing) and parametric updating steps, using a reasonable parametric model for the initialization step. Since the link function is arbitrary, except for smoothness and monotonicity constraints, we may require that estimates and parameters satisfy $\|\beta\| = 1$, $\|\hat{\beta}\| = 1$ for identifiability.

For given $\hat{\beta}, \|\hat{\beta}\| = 1$, setting $\hat{\eta}_i = \sum_{j=0}^p \hat{\beta}_j \varepsilon_j^{(i)}$, updates of the link function estimate $\hat{g}(\cdot)$ and its first derivative $\hat{g}'(\cdot)$ are obtained by smoothing (applying any reasonable scatterplot smoothing method that allows the estimation of derivatives) the scatterplot $(\hat{\eta}_i, Y_i)_{i=1,\ldots,n}$. Updates for the variance function estimate $\hat{\sigma}^2(\cdot)$ are obtained by smoothing the scatterplot $(\hat{\mu}_i, \hat{\varepsilon}_i^2)_{i=1,\ldots,n}$, where $\hat{\mu}_i = \hat{g}(\hat{\eta}_i)$ are current mean response estimates and $\hat{\varepsilon}_i^2 = (Y_i - \hat{\mu}_i)^2$ are current squared residuals. The parametric updating step then proceeds by solving the score equation (5), using the semiparametric score

(9) $$U(\beta) = \sum_{i=1}^n (Y_i - \hat{g}(\eta_i))\hat{g}'(\eta_i)\varepsilon^{(i)}/\hat{\sigma}^2(\hat{g}(\eta_i)).$$

This leads to the solutions $\hat{\beta}$, in analogy to (7). For solutions of the score equations for both scores (6) and (9), we then obtain the regression function

estimates

$$(10) \qquad \hat{\beta}(t) = \hat{\beta}_0 + \sum_{j=1}^{p} \hat{\beta}_j \rho_j(t).$$

Matrices $D$ and $\Gamma$ are modified analogously for the SPQR case, substituting appropriate estimates.

**4. Asymptotic inference.** Given an $L^2$-integrable integral kernel function $R(s,t) \colon \mathcal{T}^2 \to \mathbb{R}$, define the linear integral operator $A_R \colon L^2(dw) \to L^2(dw)$ on the Hilbert space $L^2(dw)$ for $f \in L^2(dw)$ by

$$(11) \qquad (A_R f)(t) = \int f(s) R(s,t)\, dw(s).$$

Operators $A_R$ are compact self-adjoint Hilbert–Schmidt operators if

$$\int |R(s,t)|^2\, dw(s)\, dw(t) < \infty,$$

and can then be diagonalized [Conway (1990), page 47].

Integral operators of special interest are the autocovariance operator $A_K$ of $X$ with kernel

$$(12) \qquad K(s,t) = \mathrm{cov}\,(X(s), X(t)) = \mathbf{E}(X(s) X(t))$$

and the generalized autocovariance operator $A_G$ with kernel

$$(13) \qquad G(s,t) = \mathbf{E}\left(\frac{g'(\eta)^2}{\sigma^2(\mu)} X(s) X(t)\right).$$

Hilbert–Schmidt operators $A_R$ generate a metric in $L^2$,

$$d_R^2(f,g) = \int (f(t) - g(t))(A_R(f - g))(t)\, dw(t)$$

$$= \iint (f(s) - g(s))(f(t) - g(t)) R(s,t)\, dw(s)\, dw(t)$$

for $f, g \in L^2(dw)$, and given an arbitrary orthonormal basis $\{\rho_j, j = 1, 2, \ldots\}$, the Hilbert–Schmidt kernels $R$ can be expressed as

$$R(s,t) = \sum_{k,l} r_{kl}\, \rho_k(s) \rho_l(t)$$

for suitable coefficients $\{r_{kl}, k, l = 1, 2, \ldots\}$ [Dunford and Schwartz (1963), page 1009]. Using for any given function $h \in L^2$ the notation

$$h_{\rho,j} = \int h(s) \rho_j(s)\, dw(s)$$

and denoting the normalized eigenfunctions and eigenvalues of the operator $A_R$ by $\{\rho_j^R, \lambda_j^R, j = 1, 2, \ldots\}$, the distance $d_R$ can be expressed as

$$
\begin{aligned}
(14) \quad d_R^2(f, g) &= \sum_{k,l} r_{kl}(f_{\rho,k} - g_{\rho,k})(f_{\rho,l} - g_{\rho,l}) \\
&= \sum_k \lambda_k^R (f_{\rho^R,k} - g_{\rho^R,k})^2.
\end{aligned}
$$

In the following we use the metric $d_G$, since it allows us to derive asymptotic limits under considerably simpler conditions than for the $L^2$ metric, due to its dampening effect on higher order frequencies. For the sequence of $p_n$-truncated models (1) that we are considering,

$$
d_G^2(\hat{\beta}, \beta) = \iint (\hat{\beta}(s) - \beta(s))(\hat{\beta}(t) - \beta(t)) \mathbf{E}\left(\frac{g'(\eta)^2}{\sigma^2(\mu)} X(s) X(t)\right) dw(s)\, dw(t)
$$

is approximated by $d_{G,p}^2(\hat{\beta}, \beta) = (\hat{\beta} - \beta)^T \Gamma(\hat{\beta} - \beta)$ for each $p$.

In addition to the basic assumptions in Section 2 and usual conditions on variance and link functions, we require some technical conditions which restrict the growth of $p = p_n$ and the higher-order moments of the random coefficients $\varepsilon_j$. Additional conditions are required for the semiparametric (SPQR) case where both link and variance functions are assumed unknown and are estimated nonparametrically.

(M1) The link function $g$ is monotone, invertible and has two continuous bounded derivatives with $\|g'(\cdot)\| \leq c$, $\|g''(\cdot)\| \leq c$ for a constant $c \geq 0$. The variance function $\sigma^2(\cdot)$ has a continuous bounded derivative and there exists a $\delta > 0$ such that $\sigma(\cdot) \geq \delta$.

(M2) The number of predictor terms $p_n$ in the sequence of approximating $p_n$-truncated models (1) satisfies $p_n \to \infty$ and $p_n n^{-1/4} \to 0$ as $n \to \infty$.

(M3) It holds that [see (8), where the $\xi_{kl}$ are defined]

$$
\sum_{k_1, \ldots, k_4 = 0}^{p_n} \mathbf{E}\left(\varepsilon_{k_1} \varepsilon_{k_2} \varepsilon_{k_3} \varepsilon_{k_4} \frac{g'^4(\eta)}{\sigma^4(\mu)}\right) \xi_{k_1 k_2} \xi_{k_3 k_4} = o(n/p_n^2).
$$

(M4) It holds that

$$
\begin{aligned}
\sum_{k_1, \ldots, k_8 = 0}^{p_n} &\mathbf{E}\left(\frac{g'^4(\eta)}{\sigma^4(\mu)} \varepsilon_{k_1} \varepsilon_{k_3} \varepsilon_{k_5} \varepsilon_{k_7}\right) \\
&\times \mathbf{E}\left(\frac{g'^4(\eta)}{\sigma^4(\mu)} \varepsilon_{k_2} \varepsilon_{k_4} \varepsilon_{k_6} \varepsilon_{k_8}\right) \xi_{k_1 k_2} \xi_{k_3 k_4} \xi_{k_5 k_6} \xi_{k_7 k_8} = o(n^2 p_n^2).
\end{aligned}
$$

We are now in a position to state the central asymptotic result. Given $p = p_n$, denote by $\hat{\beta} = (\hat{\beta}_0, \ldots, \hat{\beta}_p)^T$ the solution of the estimating equations (5),

(6) and by $\beta = (\beta_0, \ldots, \beta_p)^T$ the intercept $\alpha = \beta_0$ and the first $p$ coefficients of the expansion of the parameter function $\beta(t) = \sum_{j=1}^{\infty} \beta_j \rho_j(t)$ in the basis $\{\rho_j, j \geq 1\}$.

THEOREM 4.1. *If the basic assumptions and* (M1)–(M4) *are satisfied, then*

$$(15) \qquad \frac{n(\hat{\beta} - \beta)^T \Gamma_{p_n}(\hat{\beta} - \beta) - (p_n + 1)}{\sqrt{2(p_n + 1)}} \xrightarrow{d} N(0,1) \qquad as \ n \to \infty.$$

We note that the matrix $\Gamma_{p_n}$ in Theorem 4.1 may be replaced by the empirical version $\tilde{\Gamma} = \frac{1}{n}(DD^T)$; this is a consequence of (21), (22) and Lemma 7.2 below. Whenever only the "slope" parameters $\beta_1, \beta_2, \ldots$ but not the intercept parameter $\alpha = \beta_0$ are of interest, $p_n$ is replaced by $p_n - 1$ and the $(p + 1) \times (p + 1)$ matrix $\Gamma$ is replaced by the $p \times p$ submatrix of $\Gamma$ obtained by deleting the first row/column.

To study the convergence of the estimated parameter function $\hat{\beta}(\cdot)$, we use the distance $d_G$ and the representation (14) with $R \equiv G$, coupled with the expansion

$$\hat{\beta}(t) = \sum_{j=1}^{p_n} \hat{\beta}_{\rho_j^G} \rho_j^G(t)$$

of the estimated parameter function $\hat{\beta}(\cdot)$ in the basis $\{\rho_j^G, j = 1, 2, \ldots\}$, the eigenbasis of operator $A_G$ with associated eigenvalues $\lambda_j^G$. We obtain

$$d_G^2(\hat{\beta}(\cdot), \beta(\cdot)) = \iint (\hat{\beta}(s) - \beta(s)) G(s,t)(\hat{\beta}(t) - \beta(t)) \, dw(s) \, dw(t)$$

$$= \sum_{j=1}^{p} \lambda_j^G (\hat{\beta}_{\rho_j^G} - \beta_{\rho_j^G})^2 + \sum_{j=p+1}^{\infty} \lambda_j^G \beta_{\rho_j^G}^2$$

$$= (\hat{\beta}^G - \beta^G)^T \Gamma^G (\hat{\beta}^G - \beta^G) + \sum_{j=p+1}^{\infty} \lambda_j^G \beta_{\rho_j^G}^2.$$

Here

$$\hat{\beta}^G = (\hat{\beta}_{\rho_1^G}, \ldots, \hat{\beta}_{\rho_p^G})^T, \qquad \beta^G = (\beta_{\rho_1^G}, \ldots, \beta_{\rho_p^G})^T,$$

and the diagonal matrix $\Gamma^G$ is obtained by replacing in the definition of the matrix $\Gamma$ [see (8)] the $\varepsilon_j$ by $\varepsilon_j^G$ that are given by

$$\varepsilon_j^G = \frac{g'(\eta)}{\sigma(\mu)} \int X(t) \rho_j^G(t) \, dw(t),$$

with the property

$$(16) \qquad \mathbf{E}(\varepsilon_j^G \varepsilon_k^G) = \iint G(s,t) \rho_j^G(s) \rho_k^G(t) \, dw(s) \, dw(t) = \delta_{ij} \lambda_j^G.$$

These considerations lead under appropriate moment conditions to the following:

COROLLARY 4.1.    *If the parameter function $\beta(\cdot)$ has the property that*

$$(17) \qquad \sum_{j=p+1}^{\infty} \mathbf{E}(\varepsilon_j^{G2}) \left[ \int \beta(t) \rho_j^G(t) \, dw(t) \right]^2 = o\left( \frac{\sqrt{p_n}}{n} \right),$$

*then*

$$\frac{n \iint (\hat{\beta}(s) - \beta(s))(\hat{\beta}(t) - \beta(t)) G(s,t) \, dw(s) \, dw(t) - (p_n + 1)}{\sqrt{2 p_n + 1}} \xrightarrow{d} N(0,1)$$

$$as \ n \to \infty.$$

We note that property (17) relates to the rate at which higher-order oscillations, relative to the oscillations of processes $X(t)$, contribute to the $L^2$ norm of the parameter function $\beta(\cdot)$.

In case of unknown link and variance functions (SPQR), one applies scatterplot smoothing to obtain nonparametric estimates of functions and derivatives and then obtains the parameter estimates $\hat{\beta}$ as solutions of the semiparametric score equation (9). After iteration, final nonparametric estimates of the link function $\hat{g}$, its derivative $\hat{g}'$ and of the variance function $\hat{\sigma}^2$ are obtained. We implement these nonparametric curve estimators with local linear or quadratic kernel smoothers, using a bandwidth $h$ in the smoothing step. For the following result we assume these conditions:

(R1)  The regularity conditions (M1)–(M6) and (K1)–(K3) of Chiou and Müller (1998) hold uniformly for all $p_n$.

(R2)  For the bandwidths $h$ of the nonparametric function estimates for link and variance function, $h \to 0$, $\frac{nh^3}{\log n} \to \infty$ and $\|\frac{p}{\sqrt{nh^2}} \Gamma^{-1/2}\| \to 0$ as $n \to \infty$.

The following result refers to the matrix

$$(18) \qquad \hat{\Gamma} = (\hat{\gamma}_{kl})_{1 \le k, l \le p_n}, \qquad \hat{\gamma}_{k,l} = \frac{1}{n} \sum_{i=1}^{n} \left( \frac{\hat{g}'^2(\hat{\eta}_i)}{\hat{\sigma}^2(\hat{\eta}_i)} \varepsilon_{ki} \varepsilon_{li} \right).$$

COROLLARY 4.2.    *Assume* (R1) *and* (R2) *and replace the matrix* $\Gamma$ *in* (15) *by the matrix* $\hat{\Gamma}$ *from* (18). *Then* (15) *remains valid for the semiparametric quasi-likelihood (SPQR) estimates* $\hat{\beta}$ *that are obtained as solutions of the semiparametric estimating equation* (9), *substituting nonparametrically estimated link and variance functions.*

Extending the arguments used in the proofs of Theorems 1 and 2 in Chiou and Müller ([1998](#)), and assuming additional regularity conditions as described there, we find for these nonparametric function estimates,

$$\sup_t \left| \frac{\hat{g}'^2(t)}{\hat{\sigma}^2(t)} - \frac{g'^2(t)}{\sigma^2(t)} \right| = O_p\left( \frac{\log n}{nh^3} + h^2 + \frac{\sqrt{p_n}}{h^2} \|\hat{\beta} - \beta\| \right).$$

Assuming that $h \to 0, \frac{nh^3}{\log n} \to \infty$ and $\|\frac{p}{\sqrt{n}h^2}\Gamma^{-1/2}\| \to 0$, we obtain from the boundedness of the design density of the linear predictors away from 0 and $\infty$ that

$$\frac{\hat{g}'^2(\hat{\eta})}{\hat{\sigma}^2(\hat{\eta})} = \frac{g'^2(\eta)}{\sigma^2(\eta)} + o_p(1),$$

where the $o_p$-terms are uniform in $p$ following (M2). Therefore, the matrix $\hat{\Gamma}$ approximates the elements of the matrix

$$\tilde{\Gamma} = \frac{1}{n}(DD^T) = (\tilde{\gamma}_{kl})_{1 \le k,l \le p_n}, \qquad \tilde{\gamma}_{k,l} = \frac{1}{n}\sum_{i=1}^n \left( \frac{g'^2(\eta_i)}{\sigma^2(\eta_i)} \varepsilon_{ki}\varepsilon_{li} \right)$$

uniformly in $k, l$ and $p_n$. This, together with the remarks after Theorem [4.1](#), justifies the extension to the semiparametric (SPQR) case with unknown link and variance functions. This case will be included in the following, unless noted otherwise.

A common problem of inference in regression models is testing for no regression effect, that is, $H_0 : \beta \equiv$ const, which is a special case of testing for $H_0 : \beta \equiv \beta_0$ for a given regression parameter function $\beta_0$. With the representation $\beta_0(t) = \sum \beta_{0j}\rho_j(t)$, the null hypothesis becomes $H_0 : \beta_j = \beta_{0j}$, $j = 0, 1, 2, \ldots$, and $H_0$ is rejected when the test statistic in Theorem [4.1](#) exceeds the critical value $\Phi(1 - \alpha)$, for the case of a fully specified link function. Through a judicious choice of the orthonormal basis $\{\rho_j, j = 1, 2, \ldots\}$, these tests also include null hypotheses of the type $H_0 : \int \beta(t)h_j(t)\,dw(t) = \tau_j$, $j = 1, 2, \ldots$, for a sequence of linearly independent functions $h_j$; these are transformed into an orthonormal basis by Gram–Schmidt orthonormalization, whence it is easy to see that these null hypotheses translate into $H_0 : \beta_j = \tau'_j, j = 1, 2, \ldots$, for suitable $\tau'_j$ if we use the new orthonormal basis in lieu of the $\{\rho_j, j \ge 1\}$. For alternative approaches to testing in functional regression, we refer to Fan and Lin ([1998](#)).

Another application of practical interest is the construction of confidence bands for the unknown regression parameter function $\beta$. In a finite sample situation for which $p = p_n$ is given and estimates $\hat{\beta}$ for $p$-vectors $\beta$ have been determined, an asymptotic $(1 - \alpha)$ confidence region for $\beta$ according to Theorem [4.1](#) is given by $(\hat{\beta} - \beta)^T \Gamma (\hat{\beta} - \beta) \le c(\alpha)$, where $c(\alpha) = [p + 1 + \sqrt{2(p+1)}\Phi(1-\alpha)]/n$, and $\Gamma$ may be replaced by its empirical counterparts $\tilde{\Gamma}$ or $\hat{\Gamma}$. More precisely, we have the following:

Corollary 4.3. *Denote the eigenvectors/eigenvalues of the matrix* $\Gamma$ *[see* (8)*] by* $(e_1, \lambda_1), \ldots, (e_{p+1}, \lambda_{p+1})$, *and let*

$$e_k = (e_{k1}, \ldots, e_{k,p+1})^T, \qquad \omega_k(t) = \sum_{l=1}^{p+1} \rho_l(t) e_{kl}, \qquad k = 1, \ldots, p+1.$$

*Then, for large* $n$ *and* $p_n$, *an approximate* $(1 - \alpha)$ *simultaneous confidence band is given by*

$$\hat{\beta}(t) \pm \sqrt{c(\alpha) \sum_{k=1}^{p+1} \frac{\omega_k(t)^2}{\lambda_k}}. \tag{19}$$

A practical simultaneous band is obtained by substituting estimates for $\omega_k$ and $\lambda_k$ that result from empirical matrices $\tilde{\Gamma}$ or $\hat{\Gamma}$ instead of $\Gamma$.

## 5. Simulation study and model selection.

5.1. *Model order selection.* An auxiliary parameter of importance in the estimation procedure is the number $p$ of eigenfunctions that are used in fitting the function $\beta(t)$. This number has to be chosen by the statistician. An appealing method is the Akaike information criterion (AIC), due to its affinity to increasing model orders, and, in addition, we found AIC to work well in practice. We discuss here the consistency of AIC for choosing $p$ in the context of the generalized linear model with full likelihood and known link function.

Assume the linear predictor vector $\eta_p$ consists of $n$ components $\eta_{p,i} = \sum_{j=0}^{p} \varepsilon_j^i \beta_j$, $i = 1, \ldots, n$, the vector $\hat{\eta}_p$ of the components $\hat{\eta}_{p,i} = \sum_{j=0}^{p} \varepsilon_j^i \hat{\beta}_j$ and the vector $\eta$ of the components $\sum_{j=0}^{\infty} \varepsilon_j^i \beta_j$. Let $G$ be the antiderivative of the (inverse) link function $g$ so that $Y$ has the density (in canonical form) $f_Y(y) = \exp(y\eta + a(y) - G(\eta))$. In particular, $\tilde{\sigma}^2(\eta) = g'(\eta)$. The deviance is

$$\mathcal{D} = -2\ell_n(Y, \hat{\eta}_p) + 2\ell_n(Y, g^{-1}(Y)),$$

with log-likelihood

$$\ell_n(Y, \hat{\eta}_p) = \sum_{i=1}^{n} Y_i \hat{\eta}_{i,p} - \sum_{i=1}^{n} G(\hat{\eta}_{i,p}).$$

Taylor expansion yields

$$\begin{aligned}
-2\ell_n(Y, \hat{\eta}_p) = {}&-2\ell_n(Y, \eta_p) \\
&+ 2(\nabla_{\beta_p} \ell_n(Y, \hat{\eta}_p))^T (\beta_p - \hat{\beta}_p) \\
&+ (\beta_p - \hat{\beta}_p)^T \left( \frac{\partial^2}{\partial \beta_k \partial \beta_\ell} \ell_n(Y, \tilde{\eta}_p) \right) (\beta_p - \hat{\beta}_p),
\end{aligned}$$

where the second term on the right-hand side is zero, due to the score equation, and the matrix in the quadratic form is essentially $(D^T D)$. It follows from the proof of Theorem 4.1 that the quadratic form $n(\beta_p - \hat{\beta}_p)^T (\frac{D^T D}{n})(\beta_p - \hat{\beta}_p)$ has asymptotic expectation $p$. Since

$$-2\ell_n(Y, \eta_p) = -2\ell_n(Y, \eta) - 2\sum_{i=1}^n (Y_i - g(\eta_i))(\eta_{i,p} - \eta_i)$$

$$+ \sum_{i=1}^n g'(\eta_i)(\eta_{i,p} - \eta_i)^2,$$

we arrive at

$$E(\mathcal{D}) = n \sum_{k,l=p+1} E(g'(\eta)\varepsilon_k\varepsilon_l)\beta_k\beta_l - p(1 + o(1)) + E_n$$

$$= n \sum_{k,l=p+1} E\left(\frac{g'^2(\eta)}{\hat{\sigma}^2(\eta)}\varepsilon_k\varepsilon_l\right)\beta_k\beta_l - p(1 + o(1)) + E_n,$$

where $E_n$ is an expression that does not depend on $p$.

Applying the law of large numbers, and similar considerations as in Section 7, we find $\mathcal{D}/E(\mathcal{D}) \xrightarrow{p} 1$, as long as $p$ is chosen in $(p_0, cn^{1/4})$. Next, applying results of Section 7,

$$d(\hat{\beta}(\cdot), \beta(\cdot)) = \iint (\hat{\beta}(s) - \beta(s))G(s,t)(\hat{\beta}(t) - \beta(t)) \, dw(s) \, dw(t)$$

$$= (\hat{\beta}_p - \beta_p)^T \Gamma (\hat{\beta}_p - \beta_p) + \sum_{k,j=p+1}^{\infty} \gamma_{j,k}\beta_j\beta_k$$

$$+ 2\sum_{j=1}^p \sum_{k=p+1}^{\infty} \gamma_{j,k}(\hat{\beta}_j - \beta_j)\beta_k,$$

where $\gamma_{k,l} = E(\frac{g'^2(\eta)}{\hat{\sigma}^2(\eta)}\varepsilon_k\varepsilon_l)$. We obtain $E(d(\hat{\beta}(\cdot), \beta(\cdot))) = p/n(1 + o(1)) + \sum_{k,j=p+1}^{\infty} \gamma_{j,k}\beta_j\beta_k(1 + o(1))$.

This analysis shows that the target function $d(\hat{\beta}(\cdot), \beta(\cdot))$ to be minimized is asymptotically close to $E(\mathcal{D}/n) + 2p/n$. This suggests that we are in the situation considered by Shibata (1981) for sequences of linear models with normal residuals and by Shao (1997) for the more general case. While the closeness of the target function and AIC is suggestive, a rigorous proof that the order $p_A$ selected by AIC and the order $p_d$ that minimizes the target function satisfy $p_d/p_A \to 1$ in probability as $n \to \infty$ or a stronger consistency or efficiency result requires additional analysis that is not provided here. One difficulty is that the usual normality assumption is not satisfied as one operates in an exponential family or quasi-likelihood setting.

In practice, we implement AIC and the alternative Bayesian information criterion BIC by obtaining first the deviance or quasi-deviance $\mathcal{D}(p)$, dependent on the model order $p$. This is straightforward in the quasi-likelihood or maximum likelihood case with known link function, and requires integrating the score function to obtain the analogue of the log-likelihood in the SPQR case with unknown link function. Once the deviance is obtained, we choose the minimizing argument of

$$(20) \qquad\qquad C(p) = \mathcal{D}(p) + \mathcal{P}(p),$$

where $\mathcal{P}$ is the penalty term, chosen as $2p$ for the AIC and as $p \log n$ for the BIC.

Several alternative selectors that we studied were found to be less stable and more computer intensive in simulations. These included minimization of the leave-one-out prediction error, of the leave-one-out misclassification rate via cross-validation [Rice and Silverman (1991)], and of the relative difference between the Pearson criterion and the deviance [Chiou and Müller (1998)].

5.2. *Monte Carlo study.* Besides choosing the number $p$ of components to include, an implementation of the proposed generalized functional linear model also requires choice of a suitable orthonormal basis $\{\rho_j, j = 1, 2, \dots\}$. Essentially one has two options, using a fixed standard basis such as the Fourier basis $\rho_j \equiv \varphi_j \equiv \sqrt{2}\sin(\pi j t)$, $t \in [0, 1]$, $j \geq 1$, or, alternatively, to estimate the eigenfunctions of the covariance operator $A_K$ (11), (12) from the data, with the goal of achieving a sparse representation. We implemented this second option following an algorithm for the estimation of eigenfunctions which is described in detail in Capra and Müller (1997); see also Rice and Silverman (1991). Once the number of model components $p$ has been determined, the $i$th observed process is reduced to the $p$ predictors $\varepsilon_j^{(i)} = \int X_i(t)\rho_j(t)\,dw(t)$, $j = 1, \dots, p$. We substitute the estimated eigenfunctions for the $\rho_j$ and evaluate the integrals numerically.

Once we have reduced the infinite-dimensional model (1) to its $p$-truncated approximation (3), we are in the realm of finite-dimensional generalized linear and quasi-likelihood models. The parameters $\alpha$ and $\beta_1, \dots, \beta_p$ in the $p$-truncated generalized functional model are estimated by solving the respective score equation. We adopted the weighted iterated least squares algorithm which is described in McCullagh and Nelder (1989) for the case of a generalized linear or quasi-likelihood model with known link function, and the QLUE algorithm described in Chiou and Müller (1998) for the SPQR model with unknown link function.

The purpose of our Monte Carlo study was to compare AIC and BIC as selection criteria for the order $p$, to study the power of statistical tests for regression effect in a generalized functional regression model and, finally, to investigate the behavior of the semiparametric SPQR procedure for functional regression, in comparison to the maximum or quasi-likelihood implementation with a fully specified link function. The design was as follows: Pseudo-random processes based on the first 20 functions from the Fourier base $X(t) = \sum_{j=1}^{20} \varepsilon_j \varphi_j(t)$ were generated by using normal pseudo-random variables $\varepsilon_j \sim N(0, 1/j^2)$, $j \geq 1$. Choosing $\beta_j = 1/j$, $1 \leq j \leq 3$, $\beta_0 = 1$, $\beta_j = 0$, $j > 3$, we defined $\beta(t) = \sum_{j=1}^{20} \beta_j \varphi(t)$ and $p(X(\cdot)) = g(\beta_0 + \sum_{j=1}^{20} \beta_j \varepsilon_j)$, choosing logit link [with $g(x) = \exp(x)/(1 + \exp(x))$] and c-loglog link [with $g(x) = \exp(-\exp(-x))$]. Then we generated responses $Y(X) \sim$ Binomial$(p(X), 1)$ as pseudo-Bernoulli r.v.s with probability $p(X)$, obtaining a sample $(X_i(t), Y_i)$, $i = 1, \ldots, n$. Estimation methods included generalized functional linear modeling with logit, c-loglog and unspecified (SPQR) link functions.

In results not shown here, a first finding was that the AIC performed somewhat better than BIC overall, in line with theoretical expectations, and, therefore, we used AIC in the data applications. To demonstrate the asymptotic results, in particular, Theorem 4.1, we obtained empirical power functions for data generated and analyzed with the logit link, using the test statistic $T$ on the left-hand side of (16) to test the null hypothesis of no regression effect $H_0 : \beta_j = 0$, $j = 1, 2, \ldots$. This test was implemented as a one-sided test at the 5% level, that is, rejection was recorded whenever $|T| > \Phi^{-1}(0.95)$. The average rejection rate was determined over 500 Monte Carlo runs, for sample sizes $n = 50, 200$, as a function of $\delta$, $0 \leq \delta \leq 2$, where the underlying parameter vector was as described in the preceding paragraph, multiplied by $\delta$, and is given by $(\delta, \delta, \delta/2, \delta/3)$. The resulting power functions are shown in Figure 1 and demonstrate that sample size plays a critical role.

To demonstrate the usefulness of the SPQR approach with automatic link estimation, we calculated the means of the estimated regression parameter functions $\hat{\beta}(\cdot)$ over 50 Monte Carlo runs for the following cases: In each run, 1000 samples were generated with either the logit or c-loglog link function and the corresponding functions $\beta(\cdot)$ were estimated in three different ways: Assuming a logit link, a c-loglog link and assuming no link, using the SPQR method. The resulting mean function estimates can be seen in Figure 2. One finds that misspecification of the link function can lead to serious problems with these estimates and that the flexibility of the SPQR approach entails a clear advantage over methods where a link function must be specified a priori.

**6. Application to medfly data and classification.** It is a long-standing problem in evolution and ecology to analyze the interplay of longevity and reproduction. On one hand, longevity is a prerequisite for reproduction; on the other hand, numerous articles have been written about a "cost of reproduction," which is the concept that a high degree of reproduction inflicts a damage on the organism and shortens its lifespan [see, e.g., Partridge and Harvey (1985)]. The precise nature of this cost of reproduction remains elusive.

Studies with Mediterranean fruit flies (*Ceratitis capitata*), or medflies for short, have been of considerable interest in pursuing these questions as hundreds of flies can be reared simultaneously and their daily reproduction activity can be observed by simply counting the daily eggs laid by each individual fly, in addition to recording its lifetime [Carey et al. (1998a, b)]. For each medfly, one may thus obtain a reproductive trajectory and one can then ask the operational question whether particular features of this random curve have an impact on subsequent mortality [see Müller et al. (2001) for a parametric approach and Chiou, Müller and Wang (2003) for a functional model, where the egg-laying trajectories are viewed as response]. In the present framework we cast this as the problem to predict whether a fly is short- or long-lived after an initial period of egg-laying is observed. We adopt a functional binomial regression model where the initial egg-laying trajectory is the predictor process and the subsequent longevity status of the fly is the response. Of particular interest is the shape of the parameter function $\beta(\cdot)$, as it provides an indication as to which features of the egg-laying process are associated with the longevity of a fly.

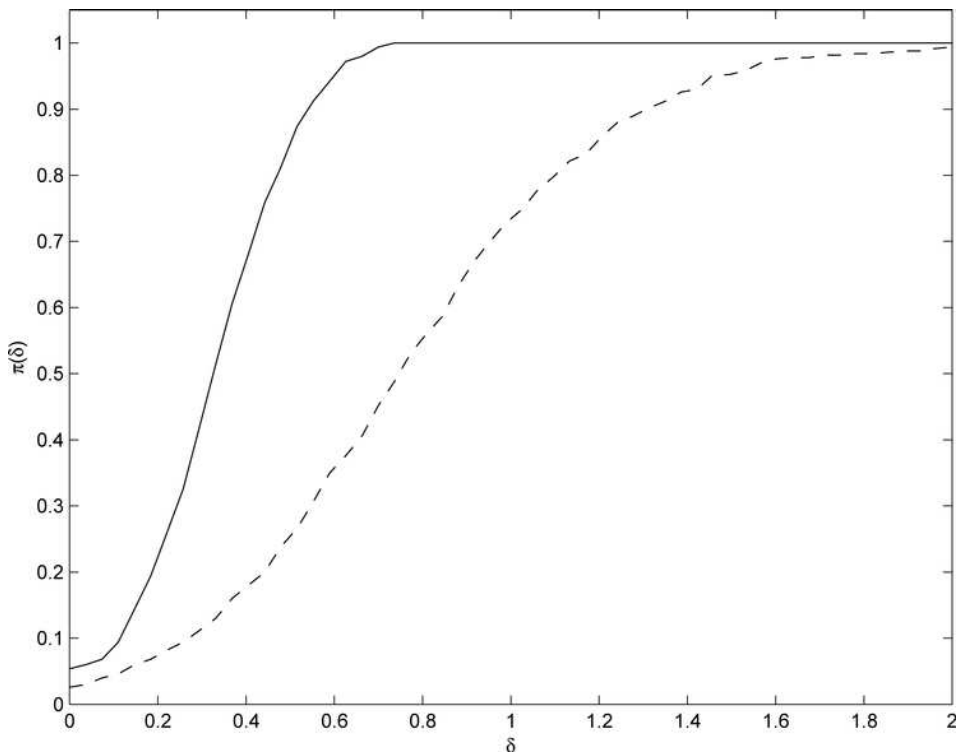

FIG. 1. *Empirical power functions for the significance test for a functional logistic regression effect at the 5% level. Based on* 500 *simulations, for sample sizes* 50 *(dashed ) and* 200 *(solid ), with* p = 3.

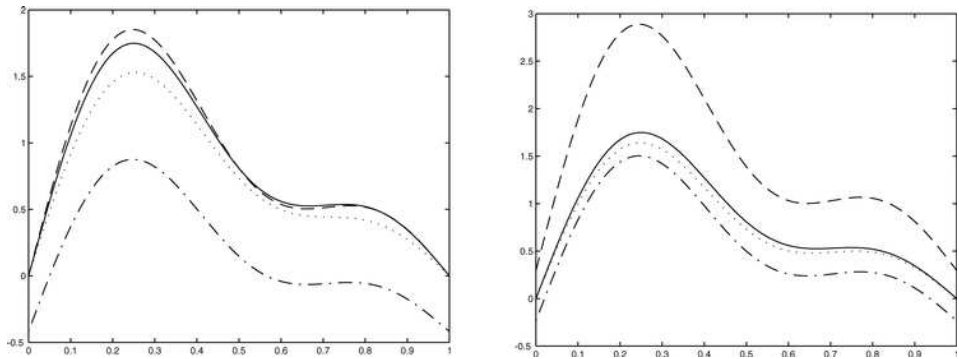

FIG. 2. *Average estimates of the regression parameter function $\beta(\cdot)$ obtained over 50 Monte Carlo runs from data generated either with the logit link (left panel) or with the c-loglog link (right panel). Each panel displays the target function (solid), and estimates obtained assuming the logit link (dashed), the c-loglog link (dash-dot) and the SPQR method incorporating nonparametric link function estimation (dotted).*

From the one thousand medflies described in Carey et al. (1998a), we select flies which lived past 34 days, providing us with a sample of 534 medflies. For prediction, we use egg-laying trajectories from 0 to 30 days, slightly smoothed to obtain the predictor processes $X_i(t)$, $t \in [0, 30]$, $i = 1, \ldots, 534$. A fly is classified as long-lived if the remaining lifetime past 30 days is 14 days or longer, otherwise as short-lived. Of the $n = 534$ flies, 256 were short-lived and 278 were long-lived. We apply the algorithm as described in the previous section, choosing the logit link, fitting a logistic functional regression.

Plotting the reproductive trajectories for the long-lived and short-lived flies separately (upper panels of Figure 3), no clear visual differences between the two groups can be discerned. Failure to visually detect differences between the two groups could result from overcrowding of these plots with too many curves, but when displaying fewer curves (lower panels of Figure 3), this remains the same. Therefore, the discrimination task at hand is difficult, as at best subtle and hard to discern differences exist between the trajectories of the two groups.

We use the Akaike information criterion (AIC) for choosing the number of model components. As can be seen from Figure 4, where the AIC criterion is shown in dependency on the model order $p$, this leads to the choice $p = 6$. The cross-validation prediction error criterion $PE = \frac{1}{n} \sum_{i=1}^{n} (Y_i - \hat{p}_i^{(-i)})^2$, where $\hat{p}_i^{(-i)}$ is the leave-one-out estimate for $p_i$, supports a similar choice. The leave-out misclassification rate estimates are, for the group of long-lived flies, 37% with logit link and 35% for the nonparametric SPQR link, while for the group of short-lived flies these are 47% for logit and 48% for SPQR, demonstrating the difficulty of classifying short-lived flies correctly.

The fitted regression parameter functions $\hat{\beta}(\cdot)$ for both logistic (logit link) and SPQR (nonparametric link) functional regression, along with simultaneous confidence bands (19), are shown in Figure 5; we find that the estimate with nonparametric link is quite close to the estimate employing the logistic link, thus providing some support for the choice of the logistic link in this case. The asymptotic confidence bands allow us to conclude that the link function has a steep rise at the right end towards age 30 days, and that the null hypothesis of no effect would be rejected.

The shape of the parameter function $\hat{\beta}(\cdot)$ highlights periods of egg-laying that are associated with increased longevity. We note that under the logit

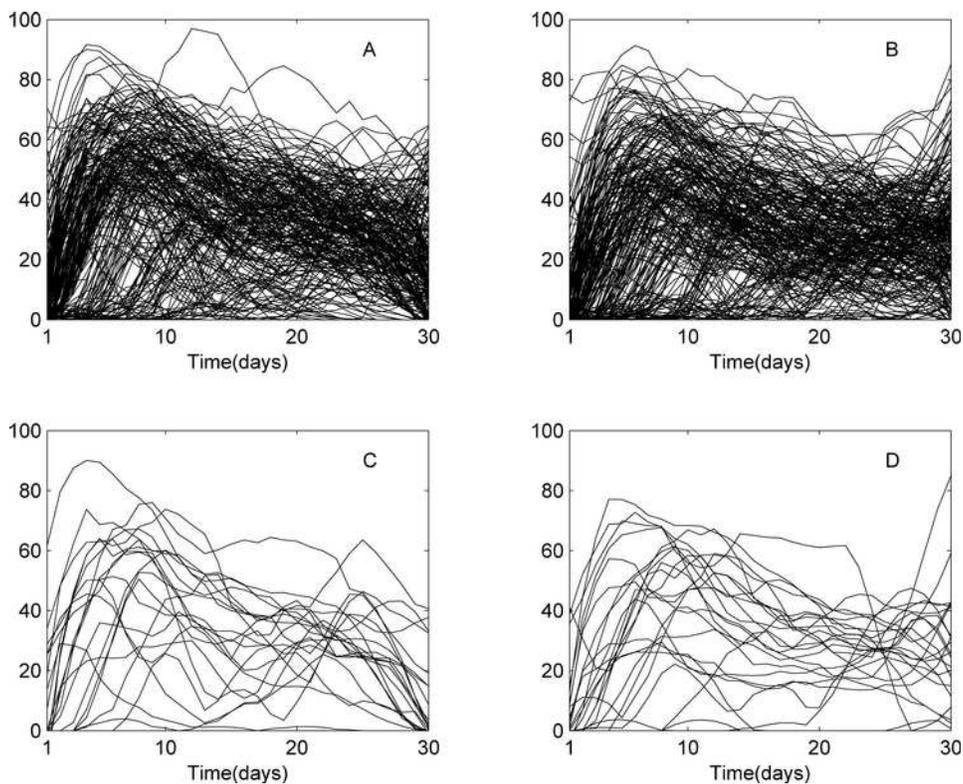

Fig. 3. *Predictor trajectories, corresponding to slightly smoothed daily egg-laying curves, for $n = 534$ medflies. The reproductive trajectories for 256 short-lived medflies are in the upper left and those for 278 long-lived medflies in the upper right panel. Randomly selected profiles from the panels above are shown in the lower panels for 50 medflies.*

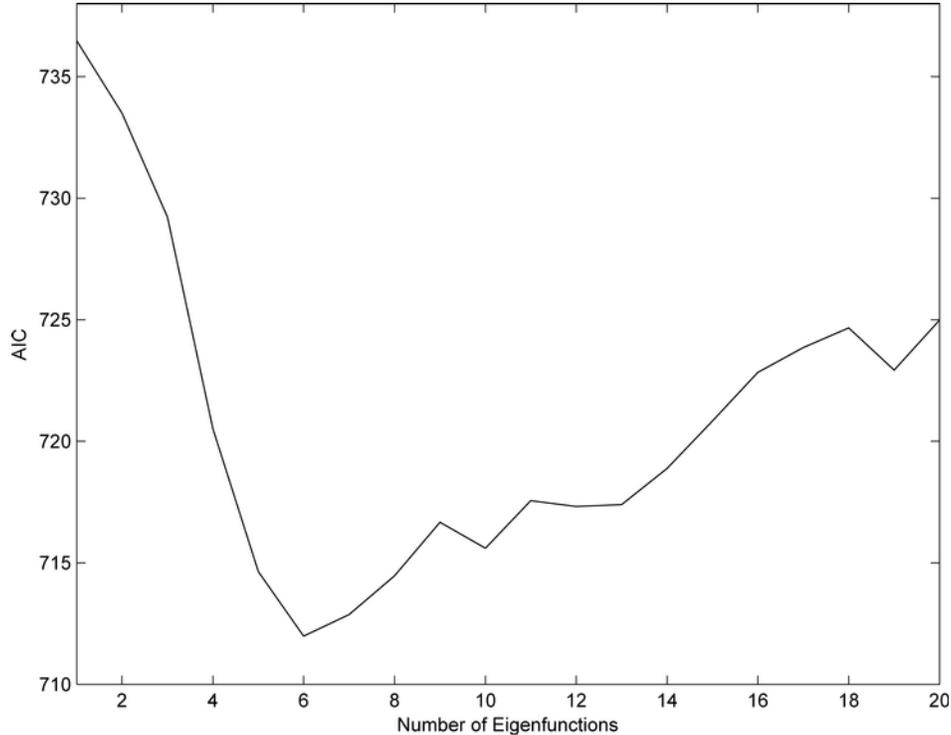

Fɪɢ. 4.  *Akaike information criterion (AIC) as a function of the number of model components p for the medfly data.*

link function, the predicted classification probability for a long-lived fly is $g(\eta) = \exp(\eta)/(1 + \exp(\eta))$. Overlaid with this expit-function, the nonparametric link function estimate that is employed in SPQR is shown in Figure 6 (choosing local linear smoothing and the bandwidth 0.55 for the smoothing steps), along with the corresponding indicator data from the last iteration step. For both links, larger linear predictors $\eta$, and therefore larger values of the parameter regression function $\beta(\cdot)$, are seen to be associated with an increased chance for longevity.

Since the parameter function is relatively large between about 12–17 days and past 26 days, we conclude that heavy reproductive activity during these periods is associated with increased longevity. In contrast, increased reproduction between 8–12 days and 20–26 days is associated with decreased longevity. A high level of late reproduction emerges as a significant and overall as the strongest indicator of longevity in our analysis. This is of biological significance since it implies that increased late reproduction is associated with increased longevity and may have a protective effect. Increased reproduction during the peak egg-laying period around 10 days has previously been associated with a cost of reproduction, an association that is supported by our analysis.

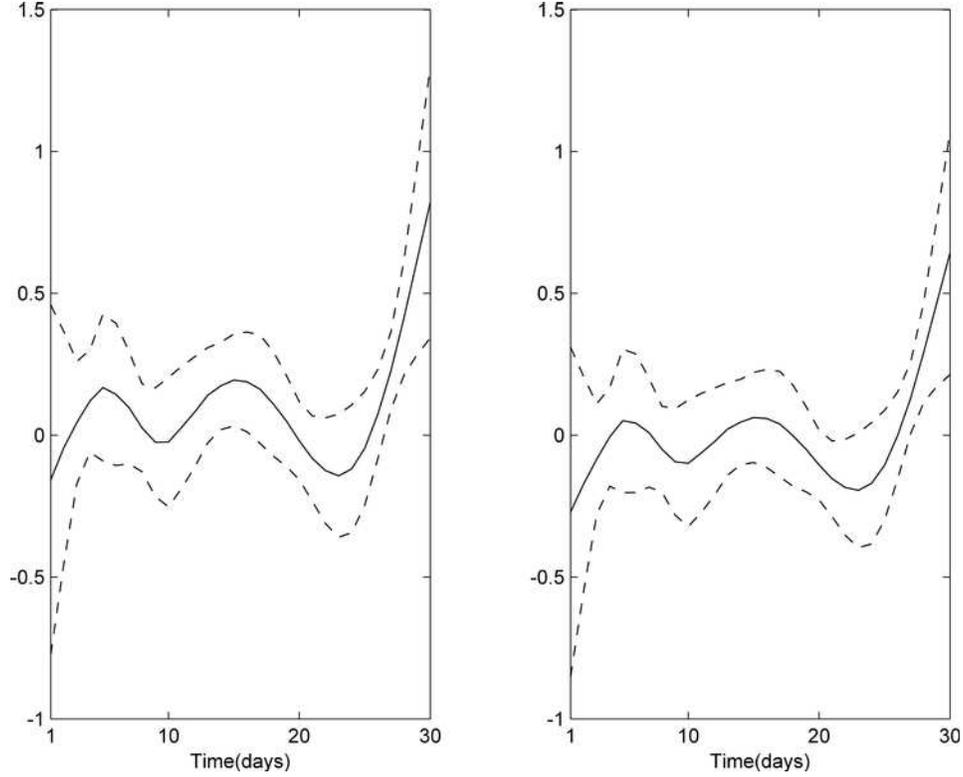

Fig. 5. *The regression parameter function estimates $\hat{\beta}(\cdot)$ (19) (solid) for the medfly classification problem, with simultaneous confidence bands (5) (dashed). Left panel: Logit link. Right panel: Nonparametric link, using the SPQR algorithm.*

**7. Proof of Theorem 4.1 and auxiliary results.** Proofs of the auxiliary results in this section are provided in the Appendix. Throughout, we assume that all assumptions of Theorem 4.1 are satisfied and work with the matrices $\Gamma = (\gamma_{kl})$, $\Xi = \Gamma^{-1} = (\xi_{kl})$, $0 \leq k, l \leq p$, defined in (8) and also with the matrix $\Xi^{1/2} =: (\xi_{kl}^{(1/2)})$, $0 \leq k, l \leq p$. We will use both versions $\sigma(\cdot)$ and $\tilde{\sigma}(\cdot)$ to represent the variance function, depending on the context, noting that $\sigma(\mu) = \tilde{\sigma}(\eta)$ and the notation $\beta$, $\hat{\beta}$ for the $(p_n + 1)$-vectors defined before Theorem 4.1 and $\beta(\cdot)$ for the parameter function.

For the first step of the proof of Theorem 4.1, we adopt the usual Taylor expansion based approach for showing asymptotic normality for an estimator which is defined through an estimating equation; see, for example, McCullagh (1983). Writing the Hessian of the quasi-likelihood as $J_\beta =$

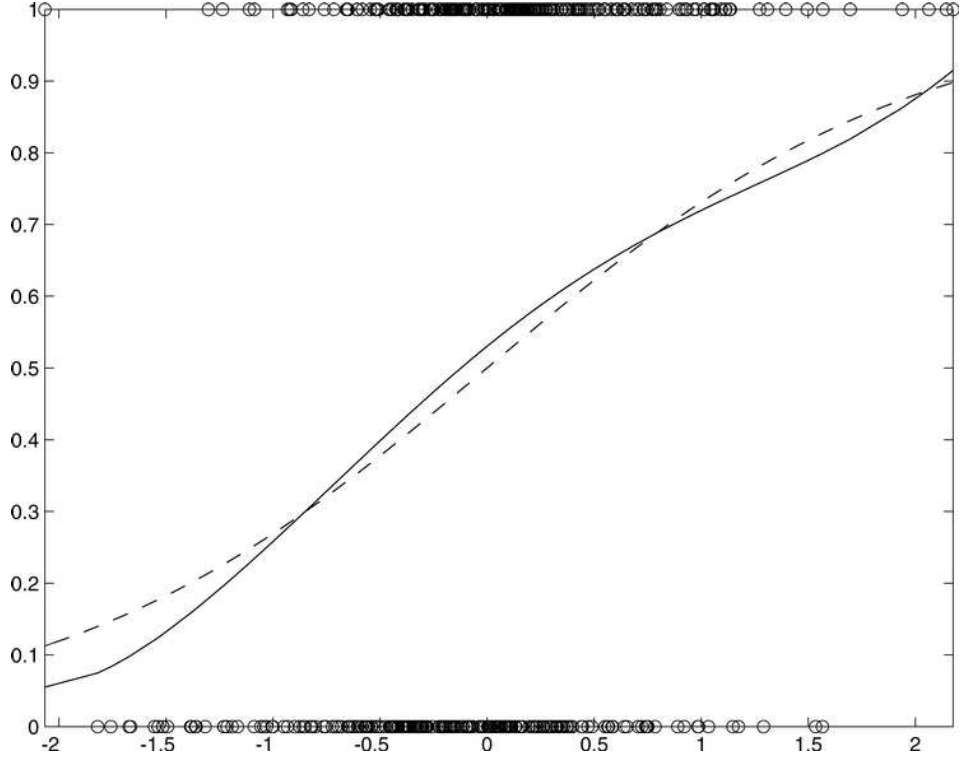

Fig. 6. *Logit link (dashed) and nonparametric link function (solid) obtained via the SPQR algorithm, with overlaid group indicators, versus level of linear predictor $\eta$.*

$\Delta_\beta U(\beta)$ and noting that

$$D^T D = \sum_{i=1}^n g'^2(\eta_i) \varepsilon^{(i)} \varepsilon^{(i)T} / \tilde{\sigma}^2(\eta_i),$$

we obtain

$$
\begin{aligned}
J_\beta &= \sum_{i=1}^n \frac{\partial}{\partial \eta_i} [g'(\eta_i) \varepsilon^{(i)} (Y_i - g(\eta_i)) / \tilde{\sigma}^2(g(\eta_i))] \cdot \Delta_\beta \eta_i \\
&= -D^T D - \sum_{i=1}^n (Y_i - g(\eta_i)) \varepsilon^{(i)} \varepsilon^{(i)T} \left\{ \frac{g''(\eta_i)}{\tilde{\sigma}^2(\eta_i)} - \frac{g'^2(\eta_i) \tilde{\sigma}^{2\prime}(\eta_i)}{\tilde{\sigma}^4(\eta_i)} \right\} \\
&= -D^T D + R, \qquad \text{say.}
\end{aligned}
$$

We aim to show that the remainder term $R$ can eventually be neglected. By a Taylor expansion, for a $\bar{\beta}$ between $\beta$ and $\hat{\beta}$,

$$
\begin{aligned}
U(\beta) &= U(\hat{\beta}) - J_{\bar{\beta}}(\hat{\beta} - \beta) = -J_{\bar{\beta}}(\hat{\beta} - \beta) \\
&= -[D^T D(\hat{\beta} - \beta) + (J_{\bar{\beta}} - J_\beta)(\hat{\beta} - \beta) + (J_\beta - D^T D)(\hat{\beta} - \beta)].
\end{aligned}
$$

Denoting the $q \times q$ identity matrix by $I_q$, this leads to

$$\sqrt{n}(\hat{\beta} - \beta) = \sqrt{n}(D^T D + (J_{\hat{\beta}} - J_{\beta}) + (J_{\beta} - D^T D))^{-1} U(\beta)$$

$$= \left( I_{p_n+1} + \left( \frac{D^T D}{n} \right)^{-1} \left( \frac{J_{\hat{\beta}} - J_{\beta}}{n} \right) + \left( \frac{D^T D}{n} \right)^{-1} \left( \frac{J_{\beta} - D^T D}{n} \right) \right)^{-1}$$

$$\times \left( \frac{D^T D}{n} \right)^{-1} \frac{U(\beta)}{\sqrt{n}}.$$

Using the matrix norm $\|M\|_2 = (\sum m_{kl}^2)^{1/2}$, we find (see Appendix for the proof) the following:

LEMMA 7.1. *As $n \to \infty$,*

$$\left\| \sqrt{n}(\hat{\beta} - \beta) - \left( \frac{D^T D}{n} \right)^{-1} \frac{U(\beta)}{\sqrt{n}} \right\|_2 = o_p(1).$$

The asymptotically prevailing term is seen to be

$$\sqrt{n}(\hat{\beta} - \beta) \sim \left( \frac{D^T D}{n} \right)^{-1} \frac{U(\beta)}{\sqrt{n}},$$

corresponding to

$$\mathcal{Z}_n = \left( \frac{D^T D}{n} \right)^{-1} \frac{D^T V^{-1/2}(Y - \mu)}{\sqrt{n}}$$

$$= \left( \frac{D^T D}{n} \right)^{-1} \frac{D^T V^{-1/2} e}{\sqrt{n}} = \left( \frac{D^T D}{n} \right)^{-1} \frac{D^T e'}{\sqrt{n}}.$$

Of interest is then the asymptotic distribution of $\mathcal{Z}_n^T \Gamma \mathcal{Z}_n$. Defining $(p+1)$-vectors $\mathcal{X}_n$ and $(p+1) \times (p+1)$-matrices $\Psi_n$ by

$$(21) \qquad \mathcal{X}_n = \frac{\Xi_n^{1/2} D^T e'}{\sqrt{n}}, \qquad \Psi_n = \Gamma_n^{1/2} \left( \frac{D^T D}{n} \right)^{-1} \Gamma_n^{1/2},$$

we may decompose this into three terms,

$$(22) \qquad \mathcal{Z}_n^T \Gamma \mathcal{Z}_n = \mathcal{X}_n^T \Psi_n^2 \mathcal{X}_n$$

$$= \mathcal{X}_n^T \mathcal{X}_n + 2\mathcal{X}_n^T (\Psi_n - I_{p_n+1}) \mathcal{X}_n$$

$$(23) \qquad \qquad + \mathcal{X}_n^T (\Psi_n - I_{p_n+1})(\Psi_n - I_{p_n+1}) \mathcal{X}_n$$

$$(24) \qquad = F_n + G_n + H_n, \qquad \text{say.}$$

The following lemma is instrumental, as it implies that in deriving the limit distribution, $G_n$ and $H_n$ are asymptotically negligible as compared to $F_n$.

Lemma 7.2. *Under the conditions*

(M3′) $p_n = o(n^{1/3})$,

(M4′) $\sum_{k_1,\ldots,k_4=0}^{p_n} \mathbf{E}(\varepsilon_{k_1}\varepsilon_{k_2}\varepsilon_{k_3}\varepsilon_{k_4}\frac{g'^4(\eta)}{\tilde\sigma^4(\eta)})\xi_{k_1 k_2}\xi_{k_3 k_4} = o(n/p_n^2)$,

*we have that*

$$\|\Psi_n - I_{p_n+1}\|_2^2 = O_p(1/p_n).$$

Note that conditions (M3′) and (M4′) are weaker than the corresponding conditions (M2) and (M3) and, therefore, will be satisfied under the basic assumptions. A consequence of Lemma 7.2 is

$$|\mathcal{X}_n^T(\Psi_n - I_{p_n+1})\mathcal{X}_n| \leq |\mathcal{X}_n\mathcal{X}_n^T|\|\Psi_n - I_{p_n+1}\|_2$$
$$= O_p(p_n)o_p(1/\sqrt{p_n}) = o_p(\sqrt{p_n}).$$

Therefore, $G_n/\sqrt{p_n} \xrightarrow{p} 0$. The bound for the term $H_n$ is completely analogous. Since we will show in Proposition 7.1 below that $(\mathcal{X}_n^T\mathcal{X}_n - (p_n + 1))/\sqrt{2p_n} \xrightarrow{d} N(0,1)$ [this implies $|\mathcal{X}_n\mathcal{X}_n^T| = O_p(p_n)$], it follows that $G_n + H_n = o_p(F_n)$ so that these terms can indeed be neglected. The proof of Theorem 4.1 will therefore be complete if we show the following:

Proposition 7.1. *As $n \to \infty$, $(\mathcal{X}_n^T\mathcal{X}_n - (p_n+1))/\sqrt{2p_n} \xrightarrow{d} N(0,1)$.*

For the proof of Proposition 7.1, we make use of

$$\mathcal{X}_n = \frac{\Xi^{1/2}D^T e'}{\sqrt{n}} = \left(\sum_{\nu=1}^n\sum_{t=0}^{p_n}\xi_{it}^{(1/2)}\frac{g'(\eta_\nu)}{\tilde\sigma(\eta_\nu)}\varepsilon_t^{(\nu)}e'_\nu/\sqrt{n}\right)_{i=0}^p$$

and

$$\sum_{k=0}^{p_n}\xi_{kt_1}^{(1/2)}\xi_{kt_2}^{(1/2)} = \xi_{t_1 t_2}$$

to obtain

$$\mathcal{X}_n^T\mathcal{X}_n = \frac{1}{n}\sum_{k=0}^{p_n}\sum_{\nu_1,\nu_2=1}^n\sum_{t_1,t_2=0}^{p_n}e'_{\nu_1}e'_{\nu_2}\frac{g'(\eta_{\nu_1})}{\tilde\sigma(\eta_{\nu_1})}\frac{g'(\eta_{\nu_2})}{\tilde\sigma(\eta_{\nu_2})}\varepsilon_{t_1}^{(\nu_1)}\varepsilon_{t_2}^{(\nu_2)}\xi_{kt_1}^{(1/2)}\xi_{kt_2}^{(1/2)}$$

$$= \frac{1}{n}\sum_{\nu=1}^n e'^2_\nu\sum_{t_1,t_2=0}^{p_n}\varepsilon_{t_1}^{(\nu)}\varepsilon_{t_2}^{(\nu)}\frac{g'^2(\eta_\nu)}{\tilde\sigma^2(\eta_\nu)}\xi_{t_1,t_2}$$

$$\quad + \frac{1}{n}\sum_{\nu_1\neq\nu_2=1}^n e'_{\nu_1}e'_{\nu_2}\frac{g'(\eta_{\nu_1})}{\tilde\sigma(\eta_{\nu_1})}\frac{g'(\eta_{\nu_2})}{\tilde\sigma(\eta_{\nu_2})}\sum_{t_1,t_2=0}^{p_n}\varepsilon_{t_1}^{(\nu_1)}\varepsilon_{t_2}^{(\nu_2)}\xi_{t_1 t_2}$$

$$= A_n + B_n, \quad \text{say.}$$

We will analyze these terms in turn and utilize the independence of the random variables associated with observations $(X_i, Y_i)$ for different values of $i$, the independence of the $e'_i$ of all $\varepsilon$'s, and $\mathbf{E}(e') = 0$, $\mathbf{E}(e'^2) = 1$.

LEMMA 7.3. *For $A_n$, it holds that*

$$\frac{A_n - (p_n + 1)}{\sqrt{p_n}} \xrightarrow{p} 0.$$

Turning now to the second term $B_n$, we show that it is asymptotically normal. Defining the r.v.s

$$W_{nj} = \sum_{k=1}^{j-1} e'_k e'_j \frac{g'(\eta_k)}{\hat{\sigma}(\eta_k)} \frac{g'(\eta_j)}{\hat{\sigma}(\eta_j)} \sum_{t_1, t_2 = 0}^{p_n} \varepsilon_{t_1}^{(k)} \varepsilon_{t_2}^{(j)} \xi_{t_1 t_2},$$

we may write

$$B_n = \frac{2}{n} \sum_{j=1}^{n} W_{nj}.$$

A key result is now the following:

LEMMA 7.4. *The random variables $\{W_{nj}, 1 \le j \le n, n \in \mathbb{N}\}$ form a triangular array of martingale difference sequences w.r.t. the filtrations $(\mathcal{F}_{nj}) = \sigma(\varepsilon_t^{(i)}, e_i, 1 \le i \le j, 0 \le t \le p_n)(1 \le j \le n, n \in \mathbb{N})$.*

Note that $\mathcal{F}_{n,j} \subset \mathcal{F}_{n+1,j}$. Lemma 7.4 implies that the r.v.s $\widetilde{W}_{nj} = \frac{2}{n\sqrt{2p_n}} W_{nj}$ also form a triangular array of martingale difference sequences. According to the central limit theorem for martingale difference sequences [Brown (1971); see also Hall and Heyde (1980), Theorem 3.2 and corollaries], sufficient conditions for the asymptotic normality $\sum_{j=1}^{n} \widetilde{W}_{nj} \xrightarrow{d} N(0, 1)$ are the conditional normalization condition and the conditional Lyapunov condition. The following two lemmas which are proved in the Appendix demonstrate that these sufficient conditions are satisfied. We note that martingale methods have also been used by Ghorai (1980) for the asymptotic distribution of an error measure for orthogonal series density estimates.

LEMMA 7.5 (Conditional normalization condition).

$$\sum_{j=1}^{n} \mathbf{E}(\widetilde{W}_{nj}^2 | \mathcal{F}_{n,j-1}) \xrightarrow{p} 1, \qquad n \to \infty.$$

Lemma 7.6 (Conditional Lyapunov condition).

$$\sum_{j=1}^{n} \mathbf{E}(\widetilde{W}_{nj}^4 | \mathcal{F}_{n,j-1}) \xrightarrow{p} 0, \qquad n \to \infty.$$

A consequence of Lemmas 7.5 and 7.6 is then $B_n/\sqrt{2 p_n} \xrightarrow{d} N(0,1)$. Together with Lemma 7.4, this implies Proposition 7.1 and, thus, Theorem 4.1.

## APPENDIX

We provide here the main arguments of the proofs of several corollaries and of the auxiliary results which were used in Section 7 for the proof of Theorem 4.1.

Proof of Corollary 4.2. Extending the arguments used in the proofs of Theorems 1 and 2 in Chiou and Müller (1998), we find for these nonparametric function estimates under (R1) that

$$\sup_t \left| \frac{\hat{g}'^2(t)}{\hat{\sigma}^2(t)} - \frac{g'^2(t)}{\sigma^2(t)} \right| = O_p\left( \frac{\log n}{nh^3} + h^2 + \frac{\sqrt{p_n}}{h^2} \|\hat{\beta} - \beta\| \right).$$

Define the matrix

$$\tilde{\Gamma} = \frac{1}{n}(DD^T) = (\tilde{\gamma}_{kl})_{1 \le k, l \le p_n}, \qquad \tilde{\gamma}_{k,l} = \frac{1}{n} \sum_{i=1}^{n} \left( \frac{g'^2(\eta_i)}{\sigma^2(\eta_i)} \varepsilon_{ki} \varepsilon_{li} \right).$$

According to (21) and (22), the result (15) remains the same when replacing $\Gamma$ by $\tilde{\Gamma}$. From (R2) and observing the boundedness of $g'^2/\sigma^2$ below and above, we obtain $\hat{\gamma}_{kl} = \tilde{\gamma}_{kl}(1 + o_p(1))$, where the $o_p$-term is uniform in $k, l$ and $p_n$. The result follows by observing that the semiparametric estimate $\hat{\beta}$ has the same asymptotic behavior as the parametric estimate, except for some minor modifications due to the identifiability constraint. □

Proof of Corollary 4.3. The asymptotic $(1 - \alpha)$ confidence ellipsoid for $\beta \in \mathbb{R}^{p+1}$ is $(\hat{\beta} - \beta)^T(\Gamma/c(\alpha))(\hat{\beta} - \beta) \le 1$. Expressing the vectors $\beta, \hat{\beta}$ in terms of the eigenvectors $e_k$ leads to the coefficients $\hat{\beta}_k^* = \sum_l e_{kl} \hat{\beta}_l$, $\beta_k^* = \sum_l e_{kl} \beta_l$, and with $\gamma_k^* = (\hat{\beta}_k^* - \beta_k^*)/\sqrt{c(\alpha)/\lambda_k}$, $\omega_k^*(t) = \omega_k(t)\sqrt{c(\alpha)/\lambda_k}$ the confidence ellipsoid corresponds to the sphere $\sum_k \gamma_k^{*2} \le 1$. To obtain the confidence band, we need to maximize $|\sum_k (\hat{\beta}_k^* - \beta_k^*)\omega_k(t)| = |\sum_k \gamma_k^* \omega_k^*(t)|$ w.r.t. $\gamma_k^*$, and subject to $\sum_k \gamma_k^{*2} \le 1$. By Cauchy–Schwarz, $|\sum_k \gamma_k^* \omega_k^*(t)| \le [\sum_k \omega_k^*(t)^2]^{1/2}$ and the maximizing $\gamma_k^*$ must be linear dependent with the vector $\omega_1^*(t), \ldots, \omega_{p+1}^*(t)$, so that the Cauchy–Schwarz inequality becomes an equality. The result then follows from the definition of the $\omega_k^*(t)$. □

Proof of Lemma 7.1.   We observe

$$\mathbf{E}\left(\left\|\frac{J_\beta - D^T D}{n}\right\|_2^2\right) = O\left(\frac{p_n^2}{n}\right) \to 0,$$

since $\|g^{(\nu)}(\cdot)\| \le c < \infty, \nu = 1, 2, \tilde{\sigma}'^2(\cdot) \le \tilde{c} < \infty$ and $\tilde{\sigma}^2(\cdot) \ge \delta > 0$ according to (M1).

Together with $p_n = o(n^{1/4})$ (M2), this implies

$$\left\|\left(\frac{D^T D}{n}\right)^{-1}\left(\frac{J_\beta - D^T D}{n}\right)\left(\frac{D^T D}{n}\right)^{-1}\frac{U(\beta)}{\sqrt{n}}\right\|_2 = o_p(1).$$

Similarly,

$$\left\|\left(\frac{D^T D}{n}\right)^{-1}\left(\frac{J_{\tilde{\beta}} - J_\beta}{n}\right)\left(\frac{D^T D}{n}\right)^{-1}\frac{U(\beta)}{\sqrt{n}}\right\| = o_p(1),$$

whence the result follows.   □

Proof of Lemma 7.2.   Note that

$$\|\Psi_n - I_{p_n+1}\|_2 \le \|\Psi_n\|_2 \|\Psi_n^{-1} - I_{p_n+1}\|_2.$$

We show that $\|\Psi_n^{-1} - I_{p_n+1}\|_2 = o_p(1)$, implying

$$\|\Psi_n\|_2 \le \|I_{p_n+1}\|_2 + \frac{\|\Psi_n^{-1} - I_{p_n+1}\|_2}{1 - \|\Psi_n^{-1} - I_{p_n+1}\|_2} \sim \|I_{p_n+1}\|_2 = \sqrt{p_n + 1}.$$

Observe that

$$\Psi_n^{-1} = \Xi_n^{1/2} \frac{1}{n} D^T D \Xi_n^{1/2}$$

$$= \left(\frac{1}{n}\sum_{j,m=0}^{p_n} \xi_{kj}^{(1/2)} \xi_{ml}^{(1/2)} \sum_{\nu=1}^n \frac{g'^2(\eta_\nu)}{\tilde{\sigma}^2(\eta_\nu)} \varepsilon_j^{(\nu)} \varepsilon_m^{(\nu)}\right)_{k,l=0}^{p_n}$$

and, therefore,

$$\mathbf{E}(\|\Psi_n^{-1} - I_{p_n+1}\|_2^2)$$

$$= \mathbf{E}\left(\sum_{k,l=0}^{p_n} \left(\frac{1}{n}\sum_{\nu_1=1}^n \sum_{j_1,m_1=0}^{p_n} \xi_{kj_1}^{(1/2)} \xi_{m_1 l}^{(1/2)} \frac{g'^2(\eta_{\nu_1})}{\sigma^2(\mu_{\nu_1})} \varepsilon_{j_1}^{(\nu_1)} \varepsilon_{m_1}^{(\nu_1)} - \delta_{kl}\right)\right.$$

$$\left. \times \left(\frac{1}{n}\sum_{\nu_2=1}^n \sum_{j_2,m_2=0}^p \xi_{kj_2}^{(1/2)} \xi_{m_2 l}^{(1/2)} \frac{g'^2(\eta_{\nu_2})}{\tilde{\sigma}^2(\eta_{\nu_2})} \varepsilon_{j_2}^{(\nu_2)} \varepsilon_{m_2}^{(\nu_2)} - \delta_{kl}\right)\right)$$

$$= O\left(\frac{p_n + 1}{n}\right) + o(1/p_n^2),$$

due to (M3$'$). Hence, by (M4$'$),

$$\|\Psi_n - I_{p_n+1}\|_2 = O_p(\sqrt{p_n})O_p\left(\frac{1}{p_n}\right) = O_p(\sqrt{1/p_n}).\qquad\square$$

PROOF OF LEMMA 7.3. Since

$$\mathbf{E}(A_n) = \frac{1}{n}\sum_{\nu=1}^{n}\sum_{t_1,t_2=0}^{p_n}\mathbf{E}(e_{\nu}'^2)\mathbf{E}\left(\varepsilon_{t_1}\varepsilon_{t_2}\frac{g'^2(\eta_\nu)}{\tilde\sigma^2(\eta_\nu)}\right)\xi_{t_1t_2} = p_n+1$$

using the definition of $\Gamma$, $\Xi = \Gamma^{-1}$ and $\mathbf{E}(e'^2) = 1$, and, similarly, by (M3),

$$\mathbf{E}(A_n^2) = o(p_n) + (p_n+1)^2 - \frac{(p_n+1)^2}{n}.$$

We find that $0 \le \mathbf{Var}(A_n) = o(p_n)$. This concludes the proof. $\square$

PROOF OF LEMMA 7.4. All random variables with upper index $j$ are independent of $\mathcal{F}_{n,j-1}$. Hence, we obtain

$$\mathbf{E}(W_{nj}|\mathcal{F}_{n,j-1}) = \sum_{i=1}^{j-1}e_i'\frac{g'(\eta_i)}{\tilde\sigma(\eta_i)}\sum_{t_1,t_2=0}^{p_n}\varepsilon_{t_1}^{(i)}\xi_{t_1t_2}\mathbf{E}\left(e_j'\frac{g'(\eta_j)}{\tilde\sigma(\eta_j)}\varepsilon_{t_2}^{(j)}|\mathcal{F}_{n,j-1}\right) = 0$$

since

$$\mathbf{E}\left(e_j'\frac{g'(\eta_j)}{\tilde\sigma(\eta_j)}\varepsilon_{t_2}^{(j)}|\mathcal{F}_{n,j-1}\right) = \mathbf{E}(e_j')\mathbf{E}\left(\frac{g'(\eta_j)}{\tilde\sigma(\eta_j)}\varepsilon_{t_2}^{(j)}\right) = 0.\qquad\square$$

PROOF OF LEMMA 7.5. We note

$$\mathbf{E}(W_{nj}^2|\mathcal{F}_{n,j-1})$$

$$= \sum_{i_1,i_2=1}^{j-1}e_{i_1}'e_{i_2}'\frac{g'(\eta_{i_1})g'(\eta_{i_2})}{\tilde\sigma(\eta_{i_1})\tilde\sigma(\eta_{i_2})}\sum_{t_1,\ldots,t_4=0}^{p_n}\varepsilon_{t_1}^{(i_1)}\varepsilon_{t_3}^{(i_2)}\xi_{t_1t_2}\xi_{t_3t_4}$$

$$\times\mathbf{E}\left(\varepsilon_{t_2}^{(j)}\varepsilon_{t_4}^{(j)}\frac{g'(\eta_j)}{\tilde\sigma^2(\eta_j)}e_j'^2|\mathcal{F}_{n,j-1}\right)$$

$$= \sum_{i_1,i_2=1}^{j-1}e_{i_1}'e_{i_2}'\frac{g'(\eta_{i_1})g'(\eta_{i_2})}{\tilde\sigma(\eta_{i_1})\tilde\sigma(\eta_{i_2})}\sum_{t_1,t_3=0}^{p_n}\varepsilon_{t_1}^{(i_1)}\varepsilon_{t_3}^{(i_2)}\xi_{t_3t_1}$$

and obtain

$$\mathbf{E}(\mathbf{E}(W_{nj}^2|\mathcal{F}_{n,j-1})) = \sum_{i=1}^{j-1}\sum_{t_1,t_3}^{p_n}\mathbf{E}\left(\frac{g'^2(\eta)}{\tilde\sigma^2(\eta)}\varepsilon_{t_1}\varepsilon_{t_3}\right)\xi_{t_3t_1}$$

$$= (j-1)(p_n+1).$$

This implies

$$\mathbf{E}\left(\sum_{j=1}^{n}\mathbf{E}(\widetilde{W}_{nj}^{2}|\mathcal{F}_{j-1})\right)\to 1,\qquad n\to\infty.$$

We are done if we can show $\mathbf{Var}(\sum_{j=1}^{n}\{\mathbf{E}(\widetilde{W}_{nj}^{2}|\mathcal{F}_{j-1})\})\to 0$. In order to obtain the second moments, we first note

$$\mathbf{E}\{\mathbf{E}(W_{nj}^{2}|\mathcal{F}_{n,j-1})\mathbf{E}(W_{nk}^{2}|\mathcal{F}_{n,k-1})\}$$

$$=\sum_{i_{1},i_{2}=1}^{j-1}\sum_{i_{3},i_{4}=1}^{k-1}\mathbf{E}\left(e_{i_{1}}'e_{i_{2}}'e_{i_{3}}'e_{i_{4}}'\frac{g'(\eta_{i_{1}})g'(\eta_{i_{2}})g'(\eta_{i_{3}})g'(\eta_{i_{4}})}{\tilde{\sigma}(\eta_{i_{1}})\tilde{\sigma}(\eta_{i_{2}})\tilde{\sigma}(\eta_{i_{3}})\tilde{\sigma}(\eta_{i_{4}})}\right)$$

$$\times\sum_{t_{1},\dots,t_{4}=0}^{p_{n}}\varepsilon_{t_{1}}^{(i_{1})}\varepsilon_{t_{2}}^{(i_{2})}\varepsilon_{t_{3}}^{(i_{3})}\varepsilon_{t_{4}}^{(i_{4})}\xi_{t_{1}t_{2}}\xi_{t_{3}t_{4}}$$

$$=\mu_{4}(k-1)\sum_{t_{1},\dots,t_{4}=0}^{p_{n}}\mathbf{E}\left(\frac{g'^{4}(\eta)}{\tilde{\sigma}^{4}(\eta)}\cdot\varepsilon_{t_{1}}\varepsilon_{t_{2}}\varepsilon_{t_{3}}\varepsilon_{t_{4}}\right)\xi_{t_{1}t_{2}}\xi_{t_{3}t_{4}}$$

$$+(j-1)(k-1)(p_{n}+1)^{2}+2(k-1)^{2}(p_{n}+1),$$

and then obtain

$$\mathbf{E}\left(\left(\sum_{j=1}^{n}\{\mathbf{E}(W_{nj}^{2}|\mathcal{F}_{n,j-1})\}\right)^{2}\right)$$

$$=\sum_{j=1}^{n}\mathbf{E}(\{\mathbf{E}(W_{nj}^{2}|\mathcal{F}_{n,j-1})\}^{2})+2\sum_{1\le k<j\le n}\mathbf{E}\{\mathbf{E}(W_{nj}^{2}|\mathcal{F}_{n,j-1})\mathbf{E}(W_{nk}^{2}|\mathcal{F}_{n,k-1})\}$$

$$=\sum_{j=1}^{n}\left[\mu_{4}(j-1)\sum_{t_{1},\dots,t_{4}=0}^{p_{n}}\mathbf{E}\left(\frac{g'^{4}(\eta)}{\tilde{\sigma}^{4}(\eta)}\cdot\varepsilon_{t_{1}}\cdots\varepsilon_{t_{4}}\right)\xi_{t_{1}t_{2}}\xi_{t_{3}t_{4}}\right.$$

$$\left.+(j-1)^{2}(p_{n}+1)^{2}+2(j-1)^{2}(p_{n}+1)\right]$$

$$+2\sum_{j=1}^{n}\sum_{k=1}^{j-1}\left((k-1)\mu_{4}\sum_{t_{1},\dots,t_{4}=0}^{p_{n}}\mathbf{E}\left(\frac{g'^{4}(\eta)}{\tilde{\sigma}^{4}(\eta)}\cdot\varepsilon_{t_{1}}\cdots\varepsilon_{t_{4}}\right)\xi_{t_{1}t_{2}}\xi_{t_{3}t_{4}}\right.$$

$$\left.+(j-1)(k-1)(p_{n}+1)^{2}+2(k-1)^{2}(p_{n}+1)\right)$$

$$= O\left( n^3 \sum_{t_1,\ldots,t_4=0}^{p_n} \mathbf{E}\left( \frac{g'^4(\eta)}{\tilde{\sigma}^4(\eta)} \cdot \varepsilon_{t_1} \cdots \varepsilon_{t_4} \right) \xi_{t_1 t_2} \xi_{t_3 t_4} \right)$$

$$+ \frac{n^4}{4}(p_n+1)^2(1+o(1)) + \frac{n^4}{6}(p_n+1)(1+o(1)).$$

Applying (M2), we infer

$$\mathbf{E}\left( \left( \sum_{j=1}^{n}\{\mathbf{E}(\widetilde{W}_{nj}^2|\mathcal{F}_{n,j-1})\} \right)^2 \right) = 1 + o(1)$$

and conclude that

$$\mathbf{Var}\left( \sum_{j=1}^{n}\{\mathbf{E}(\widetilde{W}_{nj}^2|\mathcal{F}_{n,j-1})\} \right) \to 0,$$

whence the result follows. $\square$

PROOF OF LEMMA 7.6. Combining detailed calculations of $\mathbf{E}(W_{nj}^4|\mathcal{F}_{n,j-1})$ and $\mathbf{E}(\mathbf{E}(W_{nj}^4|\mathcal{F}_{n,j-1}))$ with (M2) and (M3) leads to

$$\sum_{j=1}^{n}\mathbf{E}(\widetilde{W}_{nj}^4)$$

$$= O\left( \frac{1}{n^4 p_n^2} \right)\left[ O(n^2) \sum_{t_1,\ldots,t_8=0}^{p_n} \mathbf{E}\left( \frac{g'^4(\eta)}{\tilde{\sigma}^4(\eta)} \varepsilon_{t_1}\varepsilon_{t_3}\varepsilon_{t_5}\varepsilon_{t_7} \right) \right.$$

$$\times \mathbf{E}\left( \frac{g'^4(\eta)}{\tilde{\sigma}^4(\eta)} \varepsilon_{t_2}\varepsilon_{t_4}\varepsilon_{t_6}\varepsilon_{t_8} \right) \xi_{t_1 t_2} \xi_{t_3 t_4} \xi_{t_5 t_6} \xi_{t_7 t_8}$$

$$\left. + O(n^3) \sum_{t_3,t_4,t_7,t_8=0}^{p_n} \xi_{t_3 t_4} \xi_{t_7 t_8} \mathbf{E}\left( \frac{g'^4(\eta)}{\sigma^4(\mu)} \varepsilon_{t_3}\varepsilon_{t_4}\varepsilon_{t_7}\varepsilon_{t_8} \right) \right]$$

$$= o(1),$$

completing the proof. $\square$

**Acknowledgments.** We are grateful to an Associate Editor and two referees for helpful remarks and suggestions that led to a substantially improved version of the paper. We thank James Carey for making the medfly fecundity data available to us and Ping-Shi Wu for help with the programming.

REFERENCES

ALTER, O., BROWN, P. O. and BOTSTEIN, D. (2000). Singular value decomposition for genome-wide expression data processing and modeling. *Proc. Natl. Acad. Sci. USA* **97** 10101–10106.

ASH, R. B. and GARDNER, M. F. (1975). *Topics in Stochastic Processes.* Academic Press, New York. MR448463

BROWN, B. M. (1971). Martingale central limit theorems. *Ann. Math. Statist.* **42** 59–66. MR290428

BRUMBACK, B. A. and RICE, J. A. (1998). Smoothing spline models for the analysis of nested and crossed samples of curves (with discussion). *J. Amer. Statist. Assoc.* **93** 961–994. MR1649194

CAPRA, W. B. and MÜLLER, H.-G. (1997). An accelerated time model for response curves. *J. Amer. Statist. Assoc.* **92** 72–83. MR1436099

CARDOT, J.-M., FERRATY, F. and SARDA, P. (1999). Functional linear model. *Statist. Probab. Lett.* **45** 11–22. MR1718346

CAREY, J. R., LIEDO, P., MÜLLER, H.-G., WANG, J.-L. and CHIOU, J.-M. (1998a). Relationship of age patterns of fecundity to mortality, longevity and lifetime reproduction in a large cohort of Mediterranean fruit fly females. *J. Gerontology: Biological Sciences* **53A** B245–B251.

CAREY, J. R., LIEDO, P., MÜLLER, H.-G., WANG, J.-L. and VAUPEL, J. W. (1998b). Dual modes of aging in Mediterranean fruit fly females. *Science* **281** 996–998.

CASTRO, P. E., LAWTON, W. H. and SYLVESTRE, E. A. (1986). Principal modes of variation for processes with continuous sample curves. *Technometrics* **28** 329–337.

CHIOU, J.-M. and MÜLLER, H.-G. (1998). Quasi-likelihood regression with unknown link and variance functions. *J. Amer. Statist. Assoc.* **93** 1376–1387. MR1666634

CHIOU, J.-M. and MÜLLER, H.-G. (1999). Nonparametric quasi-likelihood. *Ann. Statist.* **27** 36–64. MR1701100

CHIOU, J.-M., MÜLLER, H.-G. and WANG, J.-L. (2003). Functional quasi-likelihood regression models with smooth random effects. *J. R. Stat. Soc. Ser. B Stat. Methodol.* **65** 405–423. MR1983755

CONWAY, J. B. (1990). *A Course in Functional Analysis*, 2nd ed. Springer, New York. MR1070713

DUNFORD, N. and SCHWARTZ, J. T. (1963). *Linear Operators. II. Spectral Theory.* Wiley, New York.

FAN, J. and LIN, S.-K. (1998). Test of significance when the data are curves. *J. Amer. Statist. Assoc.* **93** 1007–1021. MR1649196

FAN, J. and ZHANG, J.-T. (2000). Two-step estimation of functional linear models with application to longitudinal data. *J. R. Stat. Soc. Ser. B Stat. Methodol.* **62** 303–322. MR1749541

FARAWAY, J. J. (1997). Regression analysis for a functional response. *Technometrics* **39** 254–261. MR1462586

GHORAI, J. (1980). Asymptotic normality of a quadratic measure of orthogonal series type density estimate. *Ann. Inst. Statist. Math.* **32** 341–350. MR609027

HALL, P. and HEYDE, C. (1980). *Martingale Limit Theory and Its Applications.* Academic Press, New York. MR624435

HALL, P., POSKITT, D. S. and PRESNELL, B. (2001). A functional data-analytic approach to signal discrimination. *Technometrics* **43** 1–9. MR1847775

HALL, P., REIMANN, J. and RICE, J. (2000). Nonparametric estimation of a periodic function. *Biometrika* **87** 545–557. MR1789808

JAMES, G. M. (2002). Generalized linear models with functional predictors. *J. R. Stat. Soc. Ser. B Stat. Methodol.* **64** 411–432. MR1924298

MCCULLAGH, P. (1983). Quasi-likelihood functions. *Ann. Statist.* **11** 59–67. MR684863

MCCULLAGH, P. and NELDER, J. A. (1989). *Generalized Linear Models*, 2nd ed. Chapman and Hall, London. MR727836

MÜLLER, H.-G., CAREY, J. R., WU, D., LIEDO, P. and VAUPEL, J. W. (2001). Reproductive potential predicts longevity of female Mediterranean fruit flies. *Proc. R. Soc. Lond. Ser. B Biol. Sci.* **268** 445–450.

PARTRIDGE, L. and HARVEY, P. H. (1985). Costs of reproduction. *Nature* **316** 20–21.

SHAO, J. (1997). An asymptotic theory for linear model selection (with discussion). *Statist. Sinica* **7** 221–264. MR1466682

SHIBATA, R. (1981). An optimal selection of regression variables. *Biometrika* **68** 45–54. MR614940

RAMSAY, J. O. and SILVERMAN, B. W. (1997). *Functional Data Analysis.* Springer, New York.

RICE, J. A. and SILVERMAN, B. W. (1991). Estimating the mean and covariance structure nonparametrically when the data are curves. *J. Roy. Statist. Soc. Ser. B* **53** 233–243. MR1094283

STANISWALIS, J. G. and LEE, J. J. (1998). Nonparametric regression analysis of longitudinal data. *J. Amer. Statist. Assoc.* **93** 1403–1418. MR1666636

WANG, J.-L., MÜLLER, H.-G., CAPRA, W. B. and CAREY, J. R. (1994). Rates of mortality in populations of Caenorhabditis elegans. *Science* **266** 827–828.

WEDDERBURN, R. W. M. (1974). Quasi-likelihood functions, generalized linear models and the Gauss–Newton method. *Biometrika* **61** 439–447. MR375592

DEPARTMENT OF STATISTICS
UNIVERSITY OF CALIFORNIA
ONE SHIELDS AVENUE
DAVIS, CALIFORNIA 95616
USA
E-MAIL: mueller@wald.ucdavis.edu

ABT. F. ZAHLEN- U.
WAHRSCHEINLICHKEITSTHEORIE
UNIVERSITÄT ULM
89069 ULM
GERMANY